\documentclass[11pt,a4paper]{article}

\usepackage{titlesec,titletoc}
\usepackage{appendix}
\usepackage{subcaption}
\usepackage{float}
\usepackage{placeins}  
\usepackage{pdflscape}
\usepackage{longtable}
\usepackage{xcolor,colortbl,booktabs}
\usepackage{authblk}
\usepackage{dirtree}
\usepackage{algorithm}
\usepackage{algpseudocode}
\usepackage{parskip}
\usepackage{listings}
\usepackage{xcolor}
\usepackage{amsmath}
\usepackage{amssymb}
\usepackage{mathtools}
\usepackage{tcolorbox}
\usepackage{graphicx}
\usepackage{tikz}
\usetikzlibrary{intersections, calc}
\definecolor{codegreen}{rgb}{0,0.6,0}
\definecolor{codegray}{rgb}{0.5,0.5,0.5}
\definecolor{codepurple}{rgb}{0.58,0,0.82}
\definecolor{backcolour}{rgb}{1,1,1}

\lstdefinestyle{mystyle}{
    backgroundcolor=\color{backcolour},   
    commentstyle=\color{codegreen},
    keywordstyle=\color{magenta},
    numberstyle=\tiny\color{codegray},
    stringstyle=\color{codepurple},
    basicstyle=\ttfamily\footnotesize,
    breakatwhitespace=false,         
    breaklines=true,                 
    captionpos=b,                    
    keepspaces=true,                 
    numbers=left,                    
    numbersep=5pt,                  
    showspaces=false,                
    showstringspaces=false,
    showtabs=false,                  
    tabsize=2,
    frame=single
}
\makeatletter
\def\thm@space@setup{%
\thm@preskip=\parskip \thm@postskip=0pt
}
\makeatother

\tcbuselibrary{skins}
\newtcolorbox{mybox}
{
  enhanced jigsaw,
  colframe=black,
  colback=white,
  drop shadow=black!50!white,
  boxrule=0.75pt
}

\newcommand{\R}{\mathbb{R}}
\newcommand{\mb}[1]{\mathbb{#1}}
\newcommand{\mc}[1]{\mathcal{#1}}

\newcommand{\clarabel}{\textsc{clarabel}}
\newcommand{\qoco}{\textsc{qoco}}
\newcommand{\qocog}{\textsc{qocogen}}
\newcommand{\qococ}{\textsc{qoco}\textsubscript{\tiny custom}}
\newcommand{\cvxgen}{\textsc{cvxgen}}
\newcommand{\cvxpy}{\textsc{cvxpy}}
\newcommand{\cpg}{\textsc{cvxpygen}}
\newcommand{\bsocp}{\textsc{bsocp}}
\newcommand{\ecos}{\textsc{ecos}}
\newcommand{\gurobi}{\textsc{gurobi}}
\newcommand{\mosek}{\textsc{mosek}}

\usepackage{xcolor,calc}

\definecolor{BenchHighlight}{rgb}{0.6,0.8,0.6}
\newcommand{\winner}{\cellcolor{BenchHighlight}}
\usepackage[margin=0.75in]{geometry}
\usepackage{hyperref}
\hypersetup{
    colorlinks=true,
    linkcolor=blue!70,
    filecolor=PineGreen,
    urlcolor=red!50!gray!70,
    citecolor=red!50!gray!70,
    pdftitle={QOCO},
    pdfauthor={Govind M. Chari}
}

\title{\textbf{QOCO: A Quadratic Objective Conic Optimizer with Custom Solver Generation}}
\author{Govind M. Chari\! \footnote{Ph.D.\ Student, William E.\ Boeing Department of Aeronautics \& Astronautics; \texttt{gchari@uw.edu}}, \; {Beh\c{c}et}~{A\c{c}\i{}kme\c{s}e} \!\footnote{Professor, William E.\ Boeing Department of Aeronautics \& Astronautics; \texttt{behcet@uw.edu}}}  
\affil{University of Washington, Seattle, WA 98195, USA}
\date{\today}

\begin{document} 
\maketitle

\begin{abstract}
Second-order cone programs (SOCPs) with quadratic objective functions are common in optimal control and other fields. Most SOCP solvers which use interior-point methods are designed for linear objectives and convert quadratic objectives into linear ones via slack variables and extra constraints, despite the computational advantages of handling quadratic objectives directly. In applications like model-predictive control and online trajectory optimization, these SOCPs have known sparsity structures and require rapid solutions. When solving these problems, most solvers use sparse linear algebra routines, which introduce computational overhead and hinder performance. In contrast, custom linear algebra routines can exploit the known sparsity structure of problem data and be significantly faster. This work makes two key contributions: (1) the development of QOCO, an open-source C-based solver for quadratic objective SOCPs, and (2) the introduction of QOCOGEN, an open-source custom solver generator for quadratic objective SOCPs, which generates a solver written in C that leverages custom linear algebra. Both implement a primal-dual interior-point method with Mehrotra's predictor-corrector. On the benchmark problems we run, QOCO is more robust than many commonly used solvers and is faster on small- to medium-sized problems. Additionally, solvers generated by QOCOGEN are significantly faster than QOCO and are free of dynamic memory allocation making them amenable for use on embedded systems.
\end{abstract}
\section{Introduction}\label{sec:intro}

We consider the quadratic objective second-order cone program (SOCP)

\begin{equation}\label{eq:problem}
    \begin{aligned}
        \underset{x}{\text{minimize}}
        \quad & \frac{1}{2}x^\top P x + c^\top x \\
        \text{subject to}
        \quad & Gx \preceq_\mathcal{K} h \\
        \quad & Ax = b,
    \end{aligned}        
\end{equation}

with optimization variable $x \in \mathbb{R}^n$. The cost is defined by the positive semidefinite matrix $P = P^\top \succeq 0$ and vector $c \in \mathbb{R}^n$. The constraints are defined by matrices $A \in \mathbb{R}^{p \times n}$ and $G \in \mathbb{R}^{m \times n}$ and vectors $b \in \mathbb{R}^p$ and $h \in \mathbb{R}^m$. The generalized inequality $\preceq_\mathcal{K}$ denotes membership in a closed, proper cone $\mathcal{K}$, i.e. $h - Gx \in \mathcal{K}$. We restrict $\mathcal{K}$ to be the Cartesian product 

\begin{equation*}
    \mc{K} := \mc{C}_1 \times \mc{C}_2 \times \cdots \times \mc{C}_K,
\end{equation*}

where $\mc{C}_i$ is either the non-negative orthant or a second-order cone, which are self-dual. Each cone $\mc{C}_i$ corresponds to a subset of constraints, inducing a partition on $G$ and $h$, i.e., $h_i - G_ix \in \mc{C}_i$ for $i = 1,\ldots,K$. Throughout this work, we assume that Problem \ref{eq:problem} is {\bf feasible} and has a {\bf bounded} optimal objective.

Problem \ref{eq:problem} appears in various applications, including model predictive control (MPC) \cite{borelli2017mpc}, network flow optimization \cite{ford1958network}, portfolio optimization \cite{markowitz1952-ds,lobo2007portfolio}, robust optimization \cite{bental1998robust,bental1999robust}, and filter design \cite{lobo1998applications}, among others. Additionally, the solution of SOCPs is used as a subroutine in algorithms for nonconvex optimization such as sequential convex programming (SCP) \cite{malyuta2022convex, drusvyatskiy2018error,lipp2016variations,yuille2001concave}.

Problem \ref{eq:problem} can be reformulated as a SOCP with a linear objective by introducing a slack variable $t$,

\begin{equation*}
    \begin{aligned}
        \underset{x,t}{\text{minimize}}
        \quad & t + c^\top x \\
        \text{subject to}
        \quad & \left\| \begin{bmatrix}
            t - \frac{1}{2} \\
            P^{1/2} x
        \end{bmatrix}\right\|_2 \leq
        t + \frac{1}{2} \\
        \quad & Gx \preceq_\mathcal{K} h \\
        \quad & Ax = b.
    \end{aligned}        
\end{equation*}

If $P$ is large (has many rows and columns), then computing $P^{1/2}$ can be prohibitively expensive, making this reformulation impractical. Solving this reformulation can also be slow if the matrix square root $P^{1/2}$ has significantly more nonzero elements than $P$. Thus, it is advantageous for a SOCP solver to natively handle quadratic objective functions without resorting to the linear objective reformulation. 

Some SOCP solvers that can solve Problem \ref{eq:problem} directly include ${\clarabel}$ \cite{goulart2024clarabel}, $\textsc{gurobi}$ \cite{gurobi}, $\textsc{cosmo}$ \cite{garstka2021cosmo}, and $\textsc{scs}$ \cite{odonoghue2021scs}. However, $\textsc{cosmo}$ and $\textsc{scs}$ implement operator-splitting methods. These methods are preferred for large-scale problems due to their low per-iteration cost but can struggle with problem scaling, conditioning, and achieving high-accuracy solutions \cite{themelis2019supermann}. For modestly sized SOCPs, an interior-point method (IPM) is preferred, due to its robustness to scaling and conditioning of the problem, and ability to converge to high accuracy solutions \cite[Chapter~1]{Ryu2022}. One of the few open-source IPMs that can directly solve Problem \ref{eq:problem} is {\clarabel}. In addition to solving SOCPs, {\clarabel} can also handle semidefinite programs, as well as problems involving exponential cones, power cones, and generalized power cones. Additionally, it supports infeasibility detection by solving a homogeneous embedding of Problem \ref{eq:problem}. However, {\clarabel} is written in Rust. Although Rust is memory-safe due to its ownership system and can be used in embedded systems, its ecosystem is far less mature than that of C, and for legacy systems, such as aerospace software, C is still preferred.

Many engineering applications including linear MPC, nonlinear MPC \cite{rawlings1994nmpc}, portfolio backtesting, sequential quadratic programming (SQP) \cite{nocedal2006numerical,boggs1995sequential}, and SCP, require solving SOCPs with a fixed, known sparsity structures in $P, A, G$ and a constant cone $\mc{K}$. These problems often arise in real-time settings, where computational efficiency is critical. For example, sequential convex programming (SCP) iteratively solves a sequence of SOCPs with quadratic objectives and identical sparsity patterns to find stationary points of nonconvex problems. SCP has been widely applied in aerospace trajectory optimization, where real-time performance is essential. Notable applications include powered-descent guidance \cite{szmuk2020successive,kamath2023seco,chari2024fast}, in-space rendezvous \cite{berning2023chance,chari2024spacecraft}, and hypersonic entry guidance \cite{mceowen2023high,mceowen2025auto}. SCP is also the solution method used by NASA's SPLICE program for lunar landing guidance \cite{carson2019splice,Reynolds2020,kamath2023dqg,mendeck2023space}.

Most general-purpose SOCP solvers, such as ${\clarabel}$ \cite{goulart2024clarabel}, $\textsc{cosmo}$ \cite{garstka2021cosmo}, $\textsc{ecos}$ \cite{domahidi2013ecos}, $\textsc{gurobi}$ \cite{gurobi}, $\textsc{mosek}$ \cite{mosek}, and $\textsc{scs}$ rely on sparse linear algebra routines. Although sparse linear algebra allows the aforementioned solvers to solve SOCPs regardless of their sparsity structures, it can be quite slow due to the extra overhead of determining where the nonzero elements are before performing floating point operations, which leads to more CPU instructions being issued, more memory accesses, and more cache misses. However, when the sparsity structure of the problem matrices is known beforehand, it is possible to generate a custom solver that implements customized linear algebra routines tailored to the sparsity structure of the problem. Specifically, in sparse linear algebra matrices are stored as three arrays: a row/column index array, column/row pointer array, and data array. Before performing a floating point operation on a nonzero element in the matrix, sparse linear algebra routines must first index into the pointer and index arrays to determine the location of the element in the data array. In customized linear algebra, the index and pointer arrays do not exist, and we directly index into the data array, eliminating the overhead of sparse linear algebra, resulting in significantly faster operations. We discuss this further in Section \ref{subsec:custom-linalg}. 

A \textit{custom solver generator} is a tool that takes the sparsity structure of an optimization problem as an input and outputs a solver (typically in C) optimized for that specific sparsity structure without relying on sparse linear algebra. Two notable custom solver generators are {\cvxgen} \cite{mattingley2012cvxgen} and {\bsocp} \cite{dueri2014bsocp,dueri2017bsocp}. {\cvxgen}, which is academically licensed but not open-source, solves quadratic programs (QPs) rather than SOCPs (i.e. $\mc{K}$ in Problem \ref{eq:problem} is the non-negative orthant) and {\bsocp}, which is neither academically licensed nor open source, solves linear objective SOCPs (i.e. $P=0$ in Problem \ref{eq:problem}). Both have demonstrated substantial speed improvements over general-purpose solvers and have had a significant impact on aerospace trajectory optimization. For instance, after the development of lossless convexification \cite{acikmese2007lcvx, acikmese2013lossless}, {\bsocp} was flight-tested on Masten's Xombie vehicle to validate G-FOLD, a powered-descent guidance algorithm based on lossless convexification \cite{Scharf2017}. Later, {\cvxgen} was used by SpaceX for landing their Falcon 9 boosters \cite{blackmore2016autonomous}. 

A seemingly related but distinct tool is {\cpg} \cite{schaller2022embedded}. It leverages {\cvxpy} to parse the optimization problem to the standard form for solvers and generates C code for updating problem data and retrieving solutions. {\cpg} then wraps a backend solver without (or with very little) modification to its source code. Because it customizes only the parsing and solution retrieval while relying entirely on an existing solver, we do not classify {\cpg} as a custom solver generator.

\subsection{Contribution}
This work makes two key contributions: (1) {\qoco}, an open-source C-based solver for quadratic objective SOCPs, and (2) {\qocog}, an open-source custom solver generator for quadratic objective SOCPs which generates a custom solver written in C \footnote{Both are available on GitHub at \url{https://github.com/qoco-org}}. {\qoco} and solvers generated by {\qocog} implement the same algorithm as {\clarabel}, a primal-dual interior point method with Mehrotra's predictor-corrector, but directly solve Problem \ref{eq:problem} rather than a homogeneous embedding.

We demonstrate that {\qoco} is more robust than most commercial and open-source solvers, and is faster on small- to medium-sized problems. We also show that {\qococ} (the name of a solver generated by {\qocog}) is significantly faster than {\qoco}. Both solvers are easy to use since {\qoco} and {\qocog} can be called from {\cvxpy} \cite{diamond2016cvxpy,agrawal2018rewriting} and {\cpg} \cite{schaller2022embedded} respectively, allowing users to formulate optimization problems in a natural way following from math rather than manually converting the problem to the solver-required standard form. Additionally, {\qoco} can be called from C/C++, Matlab, and Python, and {\qocog} can be called from Python. Both implement a primal-dual interior point method with Mehrotra's predictor-corrector \cite{mehrotra1992pdipm}. However, they both use a variety of numerical enhancements to quickly and robustly solve the linear system that arises in the step direction computation. Specifically, we use an $LDL^\top$ factorization along with the Approximate Minimum Degree (AMD) \cite{Amestoy1996,davis2005algorithm} heuristic to permute the coefficient matrix, minimizing fill-in for the factor $L$. We also apply static and dynamic regularization to the coefficient matrix to ensure that the matrix is invertible and the factorization succeeds. After solving the linear system, we apply iterative refinement to ensure that the original, unregularized system was solved to high accuracy, further enhancing numerical stability and robustness.

\subsection{Outline}
Section \ref{sec:pdipm} introduces the primal-dual interior point method implemented in {\qoco} and {\qocog}. Section \ref{sec:implementation} explains the techniques used to ensure a stable factorization of the KKT matrix, which is essential to compute the step direction. Section \ref{sec:code-generation} discusses the advantages of custom linear algebra over sparse linear algebra, and how {\qocog} generates custom solvers which exploit the known sparsity structure of problem data. Finally, Section \ref{sec:numerical} presents extensive numerical results comparing {\qoco} and {\qococ} to existing solvers.

\subsection{Notation}
We use $(a,b)$ to denote the concatenation of vectors $a \in \mathbb{R}^a$ and $b \in \mathbb{R}^b$, resulting in a vector in $\mathbb{R}^{a+b}$. The non-negative orthant in $\mathbb{R}^l$ is denoted by $\mathbb{R}^l_+ = \{u \in \mathbb{R}^l \mid u_i \geq 0\}$, and for any $u \in \mathbb{R}^l_+$, we use $u_i$ to refer to its $i^\text{th}$ element. The $q$-dimensional second-order cone is defined as $\mathcal{Q}^{q} = \{(u_0, u_1) \in \mathbb{R} \times \mathbb{R}^{q-1} \mid \|u_1\|_2 \leq u_0\}$, where any $u \in \mathcal{Q}^q$ is written as $u = (u_0, u_1)$ with $u_0 \in \mathbb{R}$ and $u_1 \in \mathbb{R}^{q-1}$. Finally, when a vector $x$ lies in the Cartesian product of cones, i.e., $x \in \mathcal{C}_1 \times \mathcal{C}_2 \times \cdots \times \mathcal{C}_K$, we denote $x_i$ as its $i^\text{th}$ subvector belonging to cone $\mathcal{C}_i$.
\section{Primal-dual interior point method}\label{sec:pdipm}

In this section, we describe the primal-dual interior point algorithm we implement in {\qoco} and {\qocog}, as outlined in Algorithm \ref{alg:pdipm}. This algorithm is equivalent to the \texttt{coneqp} algorithm outlined in \cite{vandenberghe2010cvxopt} and follows a derivation similar to those presented in \cite[Chapter~6]{BenTal2001}, \cite{alizadeh2003second,wright1997pdipm}.

If strong duality holds for Problem \ref{eq:problem}, the optimal primal-dual solution $(x^*, s^*, y^*, z^*)$ satisfies the Karush-Kuhn-Tucker (KKT) conditions \cite[Chapter~5]{Boyd2004}, which can be written as

\begin{subequations}
    \begin{align}
        Px + c + A^\top y + G^\top z &= 0 \label{eq:kkt-eq1}\\
        Ax &= b \label{eq:kkt-eq2} \\
        Gx + s &= h \label{eq:kkt-eq3} \\
        s_i^\top z_i &= 0 \text{ for } i = 1, \ldots K \label{eq:kkt-eq4}\\
        (s, z) & \in \mc{K} \times \mc{K}. \label{eq:kkt-eq5}
    \end{align}        
\end{subequations}

The primal-dual interior point method applies a modified Newton's method to Equations \eqref{eq:kkt-eq1} - \eqref{eq:kkt-eq4}, and a line search to satisfy Equation \eqref{eq:kkt-eq5}. The modifications to Newton's method correct for linearization errors in the Newton step and bias the search directions towards the interior of $\mc{K}$. This prevents the iterates from prematurely approaching the boundary of the cone $\mc{K}$, enabling longer steps without violating Equation \eqref{eq:kkt-eq5}. It also ensures iterates do not converge to spurious solutions that satisfy Equations \eqref{eq:kkt-eq1} - \eqref{eq:kkt-eq4}, but not Equation \eqref{eq:kkt-eq5} \cite{wright1997pdipm}. 

\subsection{Central path derivation}

To understand how the algorithm biases the search directions towards the interior of $\mc{K}$, we first introduce the concept of the \textit{central path}.

Problem \ref{eq:problem} can be equivalently rewritten as

\begin{equation*}
    \begin{aligned}
        \underset{x}{\text{minimize}}
        \quad & \frac{1}{2}x^\top P x + c^\top x + \mc{I}_\mc{K}(h - Gx)\\
        \text{subject to}
        \quad & Ax = b,
    \end{aligned}        
\end{equation*}

where $\mc{I}_\mc{K}$, the indicator function of the cone $\mc{K}$, is defined as

\begin{equation*}
    \mc{I}_\mc{K}(u) = \begin{cases}
        0, \text{ for } u \in \mc{K}\\
        \infty, \text{ otherwise.}
        \end{cases}        
\end{equation*}

We then replace the indicator function, $\mc{I}_\mc{K}(u)$, with a smooth barrier function, $\phi_\mc{K}(u)$, which approaches $\infty$ as its argument approaches the boundary of the cone.

We use the barrier function

\begin{equation*}
    \phi_\mc{K}(u) = \sum_{i=1}^{K}\phi_i(u_i),
\end{equation*}

where $\phi_{i}$ is the barrier function for $\mc{C}_i$. The barrier function for the non-negative orthant and the second-order cone are

\begin{equation*}
    \phi_i(u) = 
    \begin{cases}
        -\sum_{j=1}^{l}\log u_j, \;\;\;\;\;\;\quad\quad \text{ for } \mc{C}_i = \R^l_+ \\
        -(1/2)\log (u_0^2 - u_1^\top u_1), \text{ for } \mc{C}_i = \mc{Q}^q,
    \end{cases}    
\end{equation*}

where $\log$ is the natural logarithm.

Their gradients are

\begin{equation}\label{eq:barrier-gradients}
    \nabla\phi_i(u) = 
    \begin{cases}
        (-1/u_1, \ldots, -1/u_l), \text{ for } \mc{C}_i = \R^l_+ \\
        -(u^\top J u)^{-1}Ju, \quad\;\;\;\; \text{ for } \mc{C}_i = \mc{Q}^q,
    \end{cases}    
\end{equation}

where 

\begin{equation*}
    J = \begin{bmatrix}
    1 & 0 \\
    0 & -I_{q-1}
    \end{bmatrix}.
\end{equation*}
    
After replacing the indicator function with the barrier function, we obtain

\begin{equation}\label{eq:barrier-problem}
    \begin{aligned}
        \underset{x}{\text{minimize}}
        \quad & \frac{1}{2}x^\top P x + c^\top x + \tau \sum_{i=1}^{K}\phi_{i}(h_i - G_ix)\\
        \text{subject to}
        \quad & Ax = b,
    \end{aligned}        
\end{equation}

where $\tau>0$ is a scalar. Denoting the optimal solution of Problem \ref{eq:barrier-problem} as $x^*(\tau)$, it can be shown that as $\tau \to 0$, $x^*(\tau) \to x^*$ \cite{BenTal2001}.

The KKT condition for Problem \ref{eq:barrier-problem} are

\begin{equation*}
    \begin{aligned}
        Px + c + A^\top y - \tau \sum_{i=1}^{K}G_i^\top \nabla \phi_i(h_i - G_ix) &= 0 \\
        Ax &= b \\
        h - Gx &\in \mc{K}.
    \end{aligned}        
\end{equation*}

If we define $s = h - Gx$ and $z_i = -\tau \nabla \phi_i (h_i - G_ix)$ for $i = 1, \ldots K$, we can rewrite the KKT conditions as

\begin{subequations}
    \begin{align}
        Px + c + A^\top y + G^\top z &= 0 \label{eq:central-path-1}\\
        Ax &= b \label{eq:central-path-2}\\
        Gx + s &= h \label{eq:central-path-3}\\
        z_i &= -\tau \nabla \phi_i(s_i) \text{ for } i = 1, \ldots K \label{eq:central-path-4}\\
        (s, z) & \in \mc{K} \times \mc{K}, \label{eq:central-path-5}
    \end{align}        
\end{subequations}

where the condition $z \in \mc{K}$ arises because if $u \in \mc{K}$, then $-\nabla\phi(u) \in \mc{K}$ \cite{vandenberghe2010cvxopt}.

It is desirable to write Equation \eqref{eq:central-path-4}, in a form where $s_i$ and $z_i$ appear symmetrically. To this end, we define the \textit{Jordan product} \cite{faraut1994analysis, alizadeh2003second, vandenberghe2010cvxopt}, a commutative and linear operation, for the non-negative orthant and second-order cone as

\begin{equation*}
    u_i \circ v_i = 
    \begin{cases}
        (u_{i1}v_{i1}, \ldots, u_{il}v_{il}), \quad\; \text{ for } \mc{C}_i = \R^l_+ \\
        (u_i^\top v_i, u_{i0} v_{i1} + v_{i0} u_{i1}), \text{ for } \mc{C}_i = \mc{Q}^q,
    \end{cases}    
\end{equation*}

and the Jordan product for $\mc{K}$ as

\begin{equation}\label{eq:jordan-product}
    u \circ v = (u_1 \circ v_1, u_2 \circ v_2, \ldots, u_K \circ v_K),
\end{equation}

where $u_i, v_i \in \mc{C}_i$.

The identity element $\boldsymbol{e}_i$ for cone $\mc{C}_i$ is defined as

\begin{equation*}
    \boldsymbol{e}_i = 
    \begin{cases}
        (1, 1, \ldots, 1), \text{ for } \mc{C}_i = \R^l_+ \\
        (1, 0, \ldots, 0), \text{ for } \mc{C}_i = \mc{Q}^q,
    \end{cases}    
\end{equation*}

and the identity element for $\mc{K}$ is 

\begin{equation}\label{eq:jordan-identity}
    \boldsymbol{e} = (\boldsymbol{e}_1, \ldots, \boldsymbol{e}_K).
\end{equation}

Taking the Jordan product with $s_i$ on both sides of Equation \eqref{eq:central-path-4} and substituting Equation \eqref{eq:barrier-gradients}, we obtain

\begin{subequations}
    \begin{align}
        Px + c + A^\top y + G^\top z &= 0\\
        Ax &= b\\
        Gx + s &= h \\
        s_i \circ z_i &= \tau \boldsymbol{e}_i \text{ for } i = 1, \ldots K \label{eq:cslack-cones}\\
        (s, z) & \in \mc{K} \times \mc{K}.
    \end{align}        
\end{subequations}

We then rewrite Equation \eqref{eq:cslack-cones} by stacking $s_i$ and $z_i$ and using Equations \eqref{eq:jordan-product} and \eqref{eq:jordan-identity} to obtain

\begin{subequations}
    \begin{align}
        Px + c + A^\top y + G^\top z &= 0 \label{eq:symm-central-path-1}\\
        Ax &= b \label{eq:symm-central-path-2}\\
        Gx + s &= h \label{eq:symm-central-path-3}\\
        s \circ z &= \tau \boldsymbol{e} \label{eq:symm-central-path-4}\\
        (s, z) & \in \mc{K} \times \mc{K}. \label{eq:symm-central-path-5}
    \end{align}        
\end{subequations}

The \textit{central path} is defined as the trajectory of points parameterized by $\tau>0$, satisfying Equations \eqref{eq:symm-central-path-1} - \eqref{eq:symm-central-path-5}. We denote the central path with the tuple $(x^*(\tau), s^*(\tau), y^*(\tau), z^*(\tau))$. From the above analysis, we see that the central path equations given by Equations \eqref{eq:symm-central-path-1} - \eqref{eq:symm-central-path-5} are equivalent to the optimality conditions for the barrier formulation of Problem \ref{eq:problem} given by Problem \ref{eq:barrier-problem}. It can be shown that as $\tau \to 0$, $(x^*(\tau), s^*(\tau), y^*(\tau), z^*(\tau)) \to (x^*, s^*, y^*, z^*)$ \cite{BenTal2001}.

The primal-dual interior-point method applies Newton's method to move towards the central path rather than directly towards points satisfying Equation \eqref{eq:kkt-eq4}. Since the central path lies within the cone, maintaining proximity to it allows the algorithm to take longer steps without violating Equation \eqref{eq:kkt-eq5} \cite{wright1997pdipm}.

\subsection{Nesterov-Todd scaling}

If Newton's method were applied directly to Equations \eqref{eq:central-path-1} - \eqref{eq:central-path-4}, the coefficient matrix of the resulting linear system would lack symmetry (see Appendix \ref{appendix:asymmetric-system}). As a result, it would be necessary to store the entire matrix, rather than only the upper or lower triangular portion. Additionally, solving this system would require matrix factorizations for non-symmetric matrices, which are generally more computationally expensive than those for symmetric matrices.

To derive a symmetric linear system, we introduce the following change of variables
\begin{equation*}
    \tilde{s} = W^{-\top}s, \; \tilde{z} = Wz,
\end{equation*}

where we choose $W$ such that the above transformation preserves cone membership and leaves the central path equations unchanged

\begin{equation}\label{eq:invariance}
    s \in \mc{K} \iff \tilde{s} \in \mc{K}, \quad z \in \mc{K} \iff \tilde{z} \in \mc{K}, \quad s \circ z = \tau \boldsymbol{e} \iff \tilde{s} \circ \tilde{z} = \tau \boldsymbol{e}.
\end{equation}

Using this transformation, the central path equations can be equivalently expressed as

\begin{subequations}\label{eq:scaled-central-path}
    \begin{align}
        Px + c + A^\top y + G^\top z &= 0 \\
        Ax &= b \\
        Gx + s &= h \\
        (W^{-\top} s) \circ (W z) &= \tau \boldsymbol{e} \\
        (s, z) & \in \mc{K} \times \mc{K}.
    \end{align}        
\end{subequations}

There are many matrices, $W$, which satisfy Equation \eqref{eq:invariance}. In particular, we use the \textit{Nesterov-Todd scaling} matrix \cite{nesterov1997self, nesterov1998primal}. This scaling matrix is determined based on the current iterates $s_k$ and $z_k$ and the unique point $w$ that satisfies

\begin{equation}\label{eq:scaling-point}
    \nabla^2 \phi_\mc{K}(w) s_k = z_k,
\end{equation}

where $\nabla^2 \phi_\mc{K}(w)$ is the Hessian of the barrier function evaluated at $w$. Since the barrier function is strictly convex, $\nabla^2 \phi_\mc{K}(w)$ is positive definite.

From this, we compute $W_k$ as

\begin{equation}\label{eq:scaling-matrix-equation}
    \nabla^2 \phi_\mc{K}(w)^{-1} = W_k^\top W_k.
\end{equation}

A key property that follows from Equations \eqref{eq:scaling-point} and \eqref{eq:scaling-matrix-equation} is 

\begin{equation} W_k^{-\top}s_k = W_kz_k. \end{equation}

This property allows us to define 

\begin{equation}\label{eq:lambda-definition}
    \lambda_k = W_k^{-\top}s_{k} = W_kz_{k}.
\end{equation}

For more details on how to compute the scaling point $w$ and the scaling matrix $W_k$, see \cite{vandenberghe2010cvxopt}. 

\subsection{Computing search directions}
Given the current iterate $(x_k, s_k, y_k, z_k)$ we define the \textit{duality measure} as 

\begin{equation} \mu_k = s_k^\top z_k / m, \end{equation}

which represents the average violation of the complementary slackness condition given by Equation \eqref{eq:kkt-eq4}. To compute a search direction for updating the current iterate, we apply Mehrotra's predictor-corrector method \cite{mehrotra1992pdipm}. This scheme computes a Newton step with three key objectives: reducing the duality measure, maintaining proximity to the central path, and correcting for the linearization error introduced during this step.

Applying Newton's method to Equation \eqref{eq:scaled-central-path} first requires a linearization of the central path equations around the current iterate $(x_{k},s_{k},y_{k},z_{k})$. This is done by substituting $(x,s,y,z)$ with $(x_{k}+\Delta x, s_{k}+\Delta s, y_{k}+\Delta y, z_{k}+\Delta z)$ and ignoring the higher order term $(W^{-\top} \Delta s) \circ (W \Delta z)$. This results in a linear system of the form

\begin{subequations}\label{eq:lin-central-path}
    \begin{align}
        P\Delta x + A^\top \Delta y + G^\top \Delta z &= -r_x \\
        A\Delta x &= -r_y \\
        G\Delta x + \Delta s &= -r_z \\
        (W_k^{-\top}s_{k}) \circ (W_k\Delta z) + (W_k^{-\top}\Delta s) \circ (W_kz_{k}) &= -r_s.
    \end{align}        
\end{subequations}

Here, $(r_x, r_y, r_z, r_s)$ are residual vectors that will be defined later.

Using Equation \eqref{eq:lambda-definition}, we can rewrite the system as

\begin{subequations}
    \begin{align}
        P\Delta x + A^\top \Delta y + G^\top \Delta z &= -r_x \\
        A\Delta x &= -r_y \\
        G\Delta x + \Delta s &= -r_z \\
        \lambda_k \circ (W_k\Delta z + W_k^{-\top}\Delta s) &= -r_s.
    \end{align}        
\end{subequations}

To obtain a symmetric system, we can eliminate $\Delta s$ and rewrite the equations in matrix form as

\begin{subequations}\label{eq:kkt-system}
    \begin{align}
    \begin{bmatrix}
        P & A^\top & G^\top \\
        A & 0 & 0 \\
        G & 0 & -W_k^\top W_k
    \end{bmatrix}
    \begin{bmatrix}
        \Delta x \\
        \Delta y \\
        \Delta z
    \end{bmatrix}
    &= 
    \begin{bmatrix}
        -r_x \\
        -r_y \\
        -r_z + W_k^\top (\lambda_k \backslash r_s)
    \end{bmatrix} \\
    \Delta s &= -r_z - G\Delta x,
    \end{align}        
\end{subequations}

where the operator $\backslash$ represents the inverse of the Jordan product $\circ$, meaning that $x \backslash (x \circ y) = y$.

To account for the higher order term we ignored when forming Equation \eqref{eq:lin-central-path}, we apply Mehrotra's predictor-corrector \cite{mehrotra1992pdipm}. This method consists of a predictor step, where a search direction is computed by taking a Newton step on Equation \eqref{eq:scaled-central-path} with $\tau = 0$, followed by a combined step. The combined step incorporates centering and correction terms that will help keep the next iterate close to the central path while compensating for the linearization error incurred by ignoring the $(W^{-\top} \Delta s) \circ (W \Delta z)$ term when linearizing Equation \eqref{eq:scaled-central-path}. Note that the predictor step is used to estimate this linearization error. Figure \ref{fig:pred-corr} illustrates Mehrotra's predictor-corrector.

\begin{figure}
    \centering
    \begin{tikzpicture}[scale=3]

        \fill[blue!10] (0,0) -- (3,0) -- (4,1) -- (3,2) -- (1,3) -- cycle;
        \draw[thick, blue] (0,0) -- (3,0) -- (4,1) -- (3,2) -- (1,3) -- cycle;

        \draw[red, thick] plot [smooth, tension=1] coordinates {(1.25,2) (2,1.5) (3,0)};
        
        \fill[black] (1.5,2) circle (1pt);
        \node[black] at (1.5,2.15) {$(\boldsymbol{x_{k}}, \boldsymbol{s_{k}}, \boldsymbol{y_{k}}, \boldsymbol{z_{k}})$};

        \fill[black] (2.5,1.25) circle (1pt);
        \node[black] at (3.1,1.125) {$(\boldsymbol{x_{k+1}}, \boldsymbol{s_{k+1}}, \boldsymbol{y_{k+1}}, \boldsymbol{z_{k+1}})$};

        \draw[->, ultra thick, codegreen] (1.5,2) -- (2.5,1.25);
        \draw[->, ultra thick, codepurple] (1.5,2) -- (2.5,1.75);

        \node[codegreen] at (2.9,1.4) {\textbf{Combined Step}};
        \node[codepurple] at (2.6,1.9) {\textbf{Predictor Step}};

        \node[blue] at (0.5,0.125) {\textbf{Feasible Set}};
        \node[red] at (2.375, 0.25) {\textbf{Central Path}};
        
    \end{tikzpicture}
    \caption{Mehrotra's predictor-corrector.}
    \label{fig:pred-corr}
\end{figure}

To compute the predictor step, also called the \textit{affine scaling direction}, $(\Delta x_a, \Delta s_a, \Delta y_a, \Delta z_a)$, we solve Equation \eqref{eq:kkt-system} with the residual vectors

\begin{subequations}\label{eq:affine-residuals}
    \begin{align}
        r_x &= P x_{k} + c + A^\top y_{k} + G^\top z_{k} \\
        r_y &= Ax_{k} - b \\
        r_z &= G x_{k} + s_{k} - h \\
        r_s &= \lambda_k \circ \lambda_k \label{eq:aff-rs}.
    \end{align}
\end{subequations}

Since taking a full step along the affine scaling direction can violate the constraint  $(s, z) \in \mc{K} \times \mc{K}$, we compute the centering parameter $\sigma \in [0, 1]$ as

\begin{equation*}
    \begin{split}
        \alpha &= \sup\{\alpha \in [0, 1] \; | \; (s_{k}, z_{k}) + \alpha (\Delta s_a, \Delta z_a) \in \mc{K} \times \mc{K}\} \\
        \rho &= \frac{(s_{k} + \alpha \Delta s_a)^\top (z_{k} + \alpha \Delta z_a)}{s_{k}^\top z_{k}} \\
        \sigma &= \max\{0, \min\{1, \rho\}^3\}.
    \end{split}
\end{equation*}

The centering parameter balances the need to reduce the duality measure while maintaining proximity to the central path, allowing for a longer step in the next iteration \cite{wright1997pdipm}. 

Finally, we compute the combined direction, which includes both predictor and corrector information. This direction, $(\Delta x, \Delta s, \Delta y, \Delta z)$ is obtained by solving Equation \eqref{eq:kkt-system} with the residual vectors

\begin{subequations}\label{eq:combined-residuals}
    \begin{align}
        r_x &= P x_{k} + c + A^\top y_{k} + G^\top z_{k} \\
        r_y &= Ax_{k} - b \\
        r_z &= G x_{k} + s_{k} - h \\
        r_s &= \lambda_k \circ \lambda_k + (W_k^{-\top} \Delta s_a) \circ (W_k \Delta z_a) - \sigma \mu_k\boldsymbol{e}
    \end{align}
\end{subequations}

This direction moves the next iterate closer to $(x^*(\sigma \mu_k), s^*(\sigma \mu_k), y^*(\sigma \mu_k), z^*(\sigma \mu_k))$, a point on the central path where the duality measure is reduced by a factor of $\sigma$. Note that the residual vector $r_s$ has two additional terms when compared to Equation \eqref{eq:aff-rs}: $(W_k^{-\top} \Delta s_a) \circ (W_k \Delta z_a)$ and $- \sigma \mu_k\boldsymbol{e}$. The former is to correct for linearization error in Newton's method, and the latter is to maintain proximity to the central path and arises when applying Newton's method to Equation \eqref{eq:symm-central-path-4} with $\tau = \sigma \mu_k$.

We then compute the next iterate by performing a line search to ensure $(s_{k+1}, z_{k+1}) \in \mc{K} \times \mc{K}$. Specifically, we update the iterate as

\begin{equation*}
    (x_{k+1},s_{k+1},y_{k+1},z_{k+1}) = (x_{k},s_{k},y_{k},z_{k}) + \alpha (\Delta x,\Delta s,\Delta y,\Delta z)
\end{equation*}

where

\begin{equation*}
    \alpha = \sup\left\{\alpha \in [0, 1] \; | \; (s_k, z_k) + \frac{\alpha}{0.99} ( \Delta s, \Delta z) \in \mc{K} \times \mc{K}\right\}
\end{equation*}

The factor of $0.99$ ensures that $(s_{k+1}, z_{k+1})$ remains strictly in the interior of $\mc{K} \times \mc{K}$.

\subsection{Initialization}
To initialize the primal-dual interior point method, we follow the approach outlined in \cite{vandenberghe2010cvxopt}. The initialization procedure consists of two steps. First, we find $(x, s, y, z)$ that satisfy Equations \eqref{eq:symm-central-path-1} - \eqref{eq:symm-central-path-3}. If $s$ and $z$ do not satisfy Equation \eqref{eq:symm-central-path-5}, a correction is applied to ensure they are in $\mc{K}$.

We begin by solving the optimization problem

\begin{equation*}
    \begin{aligned}
        \underset{x}{\text{minimize}}
        \quad & \frac{1}{2}x^\top P x + c^\top x + \|Gx - h\|_2^2\\
        \text{subject to}
        \quad & Ax = b.
    \end{aligned}        
\end{equation*}

Since this is an equality-constrained quadratic program, the KKT conditions correspond to the linear system

\begin{equation}\label{eq:initialize-system}
\begin{bmatrix}
    P & A^\top & G^\top \\
    A & 0 & 0 \\
    G & 0 & -I
\end{bmatrix} 
\begin{bmatrix}
    x \\
    y \\
    z
\end{bmatrix} = 
\begin{bmatrix}
    -c \\
    b \\
    h
\end{bmatrix}.
\end{equation}

We then set $s = -z$. This results in the tuple $(x, s, y, z)$ satisfying Equations \eqref{eq:symm-central-path-1} - \eqref{eq:symm-central-path-3}, but $s$ and $z$ may not necessarily satisfy Equation \eqref{eq:symm-central-path-5}.

We then initialize $x_0$ and $y_0$ as $x_0 = x$, $y_0 = y$, and initialize $s_0$ as 

\begin{equation*}
    s_0 = \begin{cases}
        s, \qquad\qquad\quad\;\; \text{ if } s \in \mc{K} \\
        s + (1 + \alpha_s)\boldsymbol{e} \;\; \text{ otherwise},
    \end{cases}
\end{equation*}

where $\alpha_s = \inf\{\alpha \; | \; s + \alpha\boldsymbol{e} \in \mc{K}\}$ and $\boldsymbol{e}$ is defined in Equation \eqref{eq:jordan-identity}. We initialize $z_0$ as 

\begin{equation*}
    z_0 = \begin{cases}
        z, \qquad\qquad\quad\;\; \text{ if } z \in \mc{K} \\
        z + (1 + \alpha_z)\boldsymbol{e}, \text{ otherwise}
    \end{cases}
\end{equation*}

where $\alpha_z = \inf\{\alpha \; | \; z + \alpha\boldsymbol{e} \in \mc{K}\}$. Note that $\alpha_{s}$ and  $\alpha_{z}$ are the minimum scalings of $\boldsymbol{e}$ that need to be added to $s$ and $z$ respectively to move them to the boundary of $\mc{K}$, but we scale $\boldsymbol{e}$ by $(1+\alpha_{s})$ and  $(1+\alpha_{z})$ to ensure that $s_0$ and $z_0$ are strictly in the interior of $\mc{K}$.

\begin{algorithm}[H]
    \caption{Primal-dual interior point method}
    \label{alg:pdipm}
    \begin{algorithmic}[1]
    \Ensure Optimal solution $(x^*,s^*,y^*,z^*)$

    \State \textbf{Initialization:}
    \State \quad Solve KKT system \eqref{eq:initialize-system} for $x, y, z$
    \State \quad $x_0 \gets x,\; y_0 \gets y$
    \State \quad $s \gets -z$

    \If{$s \in \mc{K}$}
        \State $s_0 \gets s$
    \Else
        \State $\alpha_s \gets \inf\{\alpha \;|\; s + \alpha\boldsymbol{e} \in \mc{K}\}$
        \State $s_0 \gets s + (1+\alpha_s)\boldsymbol{e}$
    \EndIf

    \If{$z \in \mc{K}$}
        \State $z_0 \gets z$
    \Else
        \State $\alpha_z \gets \inf\{\alpha \;|\; z + \alpha\boldsymbol{e} \in \mc{K}\}$
        \State $s_0 \gets z + (1+\alpha_z)\boldsymbol{e}$
    \EndIf

    \Repeat
        \State Compute duality measure:
        \State \quad $\mu_k \gets s_{k}^\top z_{k}/m$

        \State Compute Nesterov-Todd scaling:
        \State \quad Construct $W_k$
        \State \quad Calculate $\lambda_k = W_k^{-\top}s_k = W_kz_k$
        
        \State Compute affine direction $(\Delta x_a,\Delta s_a,\Delta y_a,\Delta z_a)$:
        \State \quad Solve KKT system \eqref{eq:kkt-system} with residuals \eqref{eq:affine-residuals}
        
        \State Calculate centering parameter:
        \State \quad $\alpha \gets \sup\{\alpha \in [0,1] \;|\; (s_k,z_k)+\alpha(\Delta s_a,\Delta z_a) \in \mc{K}\}$
        \State \quad $\rho \gets \frac{(s_k+\alpha\Delta s_a)^\top(z_k+\alpha\Delta z_a)}{s_k^\top z_k}$
        \State \quad $\sigma \gets \max\{0, \min\{1,\rho\}^3\}$
        
        \State Compute combined direction $(\Delta x,\Delta s,\Delta y,\Delta z)$:
        \State \quad Solve \eqref{eq:kkt-system} with residuals \eqref{eq:combined-residuals}
        
        \State Compute step-size:
        \State \quad $\alpha \gets \sup\left\{\alpha \in [0,1] \;|\; (s_k,z_k) + \frac{\alpha}{0.99}(\Delta s,\Delta z) \in \mc{K}\right\}$
        
        \State Update iterates:
        \State \quad $x_{k+1} \gets x_k + \alpha\Delta x$
        \State \quad $s_{k+1} \gets s_k + \alpha\Delta s$ 
        \State \quad $y_{k+1} \gets y_k + \alpha\Delta y$
        \State \quad $z_{k+1} \gets z_k + \alpha\Delta z$
        
    \Until{Stopping criteria \eqref{eq:stopping-criteria} satisfied.}
    \end{algorithmic}
\end{algorithm}
    
\section{Algorithm implementation}\label{sec:implementation}
This section describes the implementation of the primal-dual interior point method in {\qoco} and {\qocog}. We focus on three key aspects: the matrix factorization used to solve the system in \eqref{eq:kkt-system}, numerical techniques for ensuring a robust factorization and solution of \eqref{eq:kkt-system}, and the stopping criteria. These aspects apply to both {\qoco} and {\qocog}, but while both use the same underlying matrix factorization, their implementations are different. These differences will be discussed in Section \ref{subsec:custom-linalg}.

\subsection{Linear system solve}
The coefficient matrix in \eqref{eq:kkt-system}, which we will refer to as the Karush-Kuhn-Tucker (KKT) matrix is given as

\begin{equation}\label{eq:kkt-matrix}
    K = 
    \begin{bmatrix}
        P & A^\top & G^\top \\
        A & 0 & 0 \\
        G & 0 & -W_k^\top W_k
    \end{bmatrix}.
\end{equation}

Algorithm \ref{alg:pdipm} requires solving the linear system $K \xi = r$ for two different residual vectors $r$. Since $K$ remains constant for both solves, we factorize it once and then apply two backsolves.

The KKT matrix \eqref{eq:kkt-matrix} is symmetric but indefinite, meaning it has both positive and negative eigenvalues. A common approach for factoring such matrices is the Bunch-Parlett factorization \cite{bunch1971direct}, which factors $K$ as

\begin{equation}
    \Pi K \Pi^\top = LDL^\top,
\end{equation}

where $\Pi$ is a permutation matrix, $L$ is a lower triangular matrix, and $D$ is a block diagonal matrix with $1 \times 1$ and $2 \times 2$ blocks. However, $\Pi$ must be selected during runtime with knowledge of the data in $K$, as the factorization may not exist for all permutations. Additionally, $K$ may not be invertible if, for example, $A$ is not full row rank. 

\paragraph{Static regularization}
To ensure a numerically stable factorization that exists for all permutations $\Pi$ (allowing us to select $\Pi$ based solely on the sparsity pattern of $K$), we apply \textit{static regularization} \cite{osqp,goulart2024clarabel,mattingley2012cvxgen} to the KKT matrix

\begin{equation}\label{eq:reg-kkt-matrix}
    \hat{K} = \left[
        \begin{array}{c|cc}
          P & A^\top & G^\top\\
          \hline
          A & 0 & 0 \\
          G & 0 & -W_k^\top W_k
        \end{array}
        \right] + 
        \left[
        \begin{array}{c|cc}
          \epsilon_s I & 0 & 0 \\
          \hline
          0 & -\epsilon_s I & 0 \\
          0 & 0 & -\epsilon_s I \\
        \end{array}
        \right].
\end{equation}

In {\qoco} and {\qocog}, we use $\epsilon_s=10^{-8}$.

Now $\hat{K}$ is a \textit{quasidefinite} matrix, which is a special case of a symmetric indefinite matrix where the $(1,1)$ block is positive definite and the $(2,2)$ block is negative definite \cite{vanderbei1995symmetric}.

For such matrices, the $LDL^\top$ factorization

\begin{equation}\label{eq:ldl}
    \Pi K \Pi^\top = LDL^\top,
\end{equation}

exists for any permutation $\Pi$, and $D$ will be a diagonal matrix with known signs for its diagonal elements \cite{vanderbei1995symmetric}.

A judicious choice of $\Pi$ can minimize the \textit{fill-in} for the factor $L$. Fill-in is the introduction of non-zero elements in positions of a matrix where they previously did not exist. Excessive fill-in can lead to several negative consequences, including increased storage requirements (since more non-zero elements must be stored in memory) and increased computational cost (due to more floating-point operations). Finding an optimal $\Pi$ to minimize fill-in is NP-complete \cite{Yannakakis1981}, but heuristic methods such as the Approximate Minimum Degree (AMD) ordering can effectively compute near-optimal permutations \cite{Amestoy1996}. In {\qoco} and {\qocog}, we use Tim Davis' implementation of the AMD ordering \cite{davis2005algorithm} to determine $\Pi$.

After computing $\Pi$, the $LDL^\top$ factorization consists of two phases: a symbolic and numeric factorization. The symbolic factorization determines the elimination tree, a data structure which encodes dependencies between elements in the matrix during the factorization process, and the sparsity pattern of the factor $L$. It only requires the sparsity pattern of the regularized KKT matrix $\hat{K}$, which remains unchanged throughout Algorithm \ref{alg:pdipm}. Therefore, the symbolic factorization needs to be computed only once. The numeric factorization calculates the numerical values of $L$ and $D$, which depend on the values of the nonzero elements of $\hat{K}$. Since these values change at each iteration of Algorithm \ref{alg:pdipm}, the numeric factorization must be computed at every step. 

To compute the $LDL^\top$ factorization, we use a modified version of $\textsc{qdldl}$ \cite{osqp} with \textit{dynamic regularization} for {\qoco}, and a custom $LDL^\top$ routine for {\qocog} which we will discuss in Section \ref{subsec:custom-linalg}.

\paragraph{Dynamic regularization}
Although the $LDL^\top$ factorization must theoretically exist for $\hat{K}$, due to the nature of floating-point arithmetic it is possible for diagonal elements $D_{ii}$ to be rounded to zero or be very close to zero leading to divide-by-zero errors and the factorization failing. To alleviate this and factor $\hat{K}$ in a stable manner, we apply \textit{dynamic regularization} as given by Algorithm \ref{alg:dyn-reg}, which bounds the magnitude of $D_{ii}$ away from zero, by some $\epsilon_d > 0$, and corrects for its sign if necessary \cite{goulart2024clarabel,mattingley2012cvxgen,domahidi2013ecos}. In {\qoco} and {\qocog}, we use $\epsilon_d=10^{-8}$.

\begin{algorithm}[H]
    \caption{Dynamic regularization}
    \label{alg:dyn-reg}
    \begin{algorithmic}[1]
    \State Given: $\epsilon_d > 0$
    \If{$D_{ii}$ should be positive}
        \If{$D_{ii} < 0$}
            \State $D_{ii} = \epsilon_d$
        \Else
            \State $D_{ii} = D_{ii} + \epsilon_d$
        \EndIf
    \ElsIf{$D_{ii}$ should be negative}
        \If{$D_{ii} > 0$}
        \State $D_{ii} = -\epsilon_d$
        \Else
            \State $D_{ii} = D_{ii} - \epsilon_d$
        \EndIf
    \EndIf
    \end{algorithmic}
\end{algorithm}

\paragraph{Iterative refinement}
After static and dynamic regularization we end up solving the linear system $\hat{K}\xi = r$ rather than system $K\xi=r$. To get a solution to the latter system, we apply \textit{iterative refinement} \cite[Section 2.5.1]{press2007numerical}. This process corrects the solution by iteratively reducing the residual error, $\|r-K\xi\|$, ensuring higher numerical accuracy without requiring an additional factorization. It involves iteratively solving the system

\begin{equation}\label{eq:iter-ref}
    \hat{K}\Delta \xi^k = r - K \xi^k,
\end{equation}

and updating $\xi^{k+1} = \xi^k + \Delta \xi^k$, where $\xi^0$ solves the first linear system (i.e. $\hat{K}\xi^0 = r$). Note that since we have already computed the factorization of $\hat{K}$, one iteration of Equation \eqref{eq:iter-ref} only requires backsolves and no matrix factorization.
\subsection{Stopping criterion}
Algorithm \ref{alg:pdipm} terminates when three conditions are met: primal feasibility \eqref{eq:primal-feasibility}, stationarity \eqref{eq:stationarity}, and complementary slackness \eqref{eq:compslack}. We use both absolute and relative tolerances to ensure robustness across different problem scales. Our stopping criteria is met if the following conditions hold:

\begin{subequations}\label{eq:stopping-criteria}
    \begin{align}
        \left\|\begin{bmatrix}
            Ax_k - b \\
            Gx_k + s_k - h
        \end{bmatrix}
        \right\|_\infty &\leq \epsilon_{\mathrm{abs}} + \epsilon_{\mathrm{rel}} \max \left\{\|Ax_k\|_\infty, \|b\|_\infty, \|Gx_k\|_\infty, \|h\|_\infty, \|s_k\|_\infty\right\} \label{eq:primal-feasibility} \\
        \|Px_k + A^\top y_k + G^\top z_k + c\|_\infty & \leq \epsilon_{\mathrm{abs}} + \epsilon_{\mathrm{rel}} \max \left\{\|Px_k\|_\infty, \|A^\top y_k\|_\infty, \|G^\top z_k\|_\infty, \|c\|_\infty\right\} \label{eq:stationarity} \\
        |s_k^\top z_k| & \leq \epsilon_{\mathrm{abs}} + \epsilon_{\mathrm{rel}} \max \left\{1, |p_k|, |d_k|\right\}, \label{eq:compslack}
    \end{align}
\end{subequations}

where the primal and dual objectives are

\begin{equation*}
    p_k := \frac{1}{2}x_k^\top P x_k + c^\top x_k, \quad d_k := -\frac{1}{2}x_k^\top P x_k - b^\top y_k - h^\top z_k.
\end{equation*}
\section{Sparsity exploiting custom solver}\label{sec:code-generation}
This section discusses {\qocog}, which generates an implementation of Algorithm \ref{alg:pdipm} called {\qococ}. Unlike generic solvers which use sparse linear algebra, {\qococ} uses a custom $LDL^\top$ factorization and other tailored linear algebra routines, which exploit the known sparsity structure of problem data, to solve Problem \ref{eq:problem} significantly faster.

While {\qoco} can solve instances of Problem \ref{eq:problem} with any sparsity pattern for $P, A, $ and $G$ and cone $\mc{K}$, the use of custom linear algebra in {\qococ} restricts it to instances with the same sparsity structure for $P, A, $ and $G$ and optimizes over the same cone $\mc{K}$. This makes {\qococ} useful for applications that repeatedly solve optimization problems with identical sparsity structures, such as sequential convex programming \cite{malyuta2022convex}.

Figure \ref{fig:qoco-vs-qocogen} illustrates the usage of {\qoco}, {\qocog}, and {\qococ}. In the figure, the matrices $(P,A,G)$ contain two pieces of information: their sparsity patterns $\mathrm{str}(P,A,G)$ (the location of nonzero elements) and the values of the nonzero elements $\mathrm{data}(P,A,G)$. We see that {\qoco} takes in $(P,A,G)$, which contains both the sparsity patterns and nonzero elements for the matrices, and the rest of the problem data and returns the optimal solution $(x^*, s^*, y^*, z^*)$. In contrast, {\qocog} only takes in the sparsity patterns of $(P,A,G)$ and cone $\mc{K}$, and then generates {\qococ}. This instance of {\qococ} can then solve optimization problems where the sparsity patterns of $(P,A,G)$ do not change, but the nonzero values for $(P,A,G)$ and the vector data $(c,b,h)$ can change.

Although {\qoco} avoids dynamic memory allocation during the solution process, it still relies on \texttt{setup} and \texttt{cleanup} functions for dynamic memory allocation and deallocation before and after solving the problem. In contrast, {\qococ} exclusively uses static memory allocation, making it well-suited for resource-constrained embedded systems, where dynamic memory allocation can lead to non-deterministic behavior and memory fragmentation.

\begin{figure}
    \begin{mybox}
    \centering
    \begin{tikzpicture}[scale=0.8, transform shape]
    
    \draw (0,3) rectangle node {\qoco} ++(3,1);
    \draw[->] (3,3.5) -- node[right, xshift=10pt] {$(x^*, s^*, y^*, z^*)$} ++(1,0);
    \draw[->] (-1,3.5) -- node[left, xshift=-10pt] {$(P,A,G,c,b,h,\mc{K})$} ++(1,0);

    \draw (0,0) rectangle node {\qocog} ++(3,1);    
    \draw[->] (3,0.5) -- node[above] {Generates} ++(2,0);
    
    \draw[->] (-1,0.5) -- node[left, xshift=-10pt] {$\mathrm{str}(P,A,G), \mc{K}$} ++(1,0);
    \draw (5,0) rectangle node {\qococ} ++(2.5,1);
    \draw[->] (6.25,2) -- node[above, yshift=10pt] {$\mathrm{data}(P,A,G),c,b,h$} ++(0,-1);
    \draw[->] (7.5,0.5) -- node[right, xshift=10pt] {$(x^*, s^*, y^*, z^*)$} ++(1,0);
        
    \end{tikzpicture}
    \end{mybox}
    \caption{Usage of {\qoco}, {\qocog}, and {\qococ}.}
    \label{fig:qoco-vs-qocogen}
\end{figure}

\subsection{Custom linear algebra}\label{subsec:custom-linalg}
A key difference between {\qococ} and {\qoco} is that {\qococ} uses custom linear algebra routines rather than sparse linear algebra. The main motivation for custom routines is to speed up the $LDL^\top$ factorization of the regularized KKT matrix, which is the most computationally expensive step of Algorithm \ref{alg:pdipm}. Typically, the KKT matrix is factored using a sparse linear algebra routine that stores the matrix in either \textit{compressed sparse column} (CSC) or \textit{compressed sparse row} (CSR) format. However, these sparse routines involve more than just the floating-point operations required for factorization. They also incur additional overhead to locate nonzero elements and their positions in the matrix. For example, in a CSC matrix, accessing data first requires indexing into the column pointer and row index arrays. This results in more CPU instructions, increased memory accesses, and a higher likelihood of cache misses. We provide empirical evidence supporting these claims in Appendix \ref{appendix:custom-ldl-perf}.

When the sparsity structure of the KKT matrix is known \textit{a priori} it is possible to write a custom $LDL^\top$ factorization routine which only contains code to perform the necessary floating-point operations and data accesses, eliminating the additional overhead that typical sparse linear algebra routines have. This makes the custom $LDL^\top$ factorization significantly faster than a sparse $LDL^\top$ factorization. We also note that since the AMD ordering we employ to compute the permutation matrix $\Pi$ only depends on the sparsity pattern of the KKT matrix, we can determine $\Pi$ at the time of code generation, so the generated code does not include the code to compute the AMD ordering.

A concrete example of a custom linear algebra routine for matrix-vector multiplication is given by Listing \ref{lst:custom-matvec}, where the sparsity pattern of $A$ is fixed and defined in Equation \eqref{eq:toy-smpv}. Note that the custom routine only includes the necessary floating-point operations. Although Listing \ref{lst:custom-matvec} depicts Python code, the custom linear algebra routines in {\qococ} are in C.

\begin{equation}\label{eq:toy-smpv}
    \begin{bmatrix}
        y_0 \\
        y_1 \\
        y_2
    \end{bmatrix}
     = 
    \underbrace{
    \begin{bmatrix}
        a_{00} & 0 & a_{02} \\
        0 & 0 & 0 \\
        a_{20} & 0 & 0
    \end{bmatrix}}_{A}
    \begin{bmatrix}
        x_0 \\
        x_1 \\
        x_2
    \end{bmatrix}
\end{equation}

\begin{lstlisting}[language=Python, caption=Custom matrix-vector multiplication., float, label=lst:custom-matvec]
    # Adata stores the nonzeros of A in column-major format, i.e.
    # Adata[0] = a00
    # Adata[1] = a20
    # Adata[2] = a02
    def custom_matvec(y, Adata, x):
        y[0] = Adata[0]*x[0] + Adata[2]*x[2]
        y[1] = 0
        y[2] = Adata[1]*x[0]
\end{lstlisting}

A drawback of these custom routines is that the amount of code generated increases with problem size. Consequently, both the time required to generate and compile {\qococ} and the resulting binary size increase for larger problems. On processors with small caches, such as the Raspberry Pi Compute Module 4, these large binaries can degrade performance due to increased instruction cache misses. This effect is discussed further in Section \ref{sec:numerical}. 
\subsection{Code generation}
{\qocog} is written in Python and example code for generating {\qococ} for a problem family defined by Equation \eqref{eq:sample-problem} is given by Listing \ref{lst:qocogen-usage}. The generated solver {\qococ} is customized to the sparsity pattern of the matrices $P, A, G$ which are passed to the function $\texttt{qocogen.generate\_solver()}$ as sparse \texttt{scipy} matrices.

\begin{equation}\label{eq:sample-problem}
    \begin{array}{ll}
        \underset{x}{\text{minimize}}
        & x_1^2+x_2^2+x_3^2+x_4 \\
        \mbox{subject to} & x_1+x_2=1 \\
        & x_2+x_3 = 1 \\
        & x_1 \geq 0 \\
        & \sqrt{x_3^2+x_4^2} \leq x_2
    \end{array}    
\end{equation}

The file tree for {\qococ} is given by Figure \ref{fig:file-tree}. The file \texttt{qoco\_custom.c} contains the \texttt{qoco\_custom\_solve()} function which implements Algorithm \ref{alg:pdipm}, \texttt{ldl.c} contains the custom $LDL^\top$ factorization of the regularized KKT matrix, \texttt{CMakeLists.txt} is the CMake configuration file which automates the build process \cite{CMake}, \texttt{runtest.c} gives an example of calling {\qococ} with default data and times the resulting solve, 
and the remaining files contain various helper functions needed to implement Algorithm \ref{alg:pdipm}.

\begin{lstlisting}[language=Python, caption=Generating a custom solver with {\qocog}., float, label=lst:qocogen-usage]
    import qocogen
    import numpy as np
    from scipy import sparse
    
    # Define problem data
    P = sparse.diags([2, 2, 2, 0], 0).tocsc()
    
    c = np.array([0, 0, 0, 1])
    G = -sparse.identity(4).tocsc()
    h = np.zeros(4)
    A = sparse.csc_matrix([[1, 1, 0, 0], [0, 1, 1, 0]])
    b = np.array([1, 1])
    
    l = 1
    n = 4
    m = 4
    p = 2
    nsoc = 1
    q = np.array([3])
    
    # Generate custom solver.
    qocogen.generate_solver(n, m, p, P, c, A, b, G, h, l, nsoc, q)
\end{lstlisting}

\begin{figure}[ht]
    \centering
    \fbox{\small
        \begin{minipage}{0.6\textwidth}
            \dirtree{%
            .1 qoco\_custom/.
            .2 cone.c.
            .2 cone.h.
            .2 kkt.c.
            .2 kkt.h.
            .2 ldl.c.
            .2 ldl.h.
            .2 qoco\_custom.c.
            .2 qoco\_custom.h.
            .2 runtest.c.
            .2 utils.c.
            .2 utils.h.
            .2 workspace.h.
            .2 CMakeLists.txt.
            }            
        \end{minipage}
    }
    \caption{{\qococ} file structure.}
    \label{fig:file-tree}
\end{figure}

Listing \ref{lst:qococ-usage} provides an example of calling {\qococ} from C, updating the problem data, and solving it again. 

\begin{lstlisting}[language=C, caption=Calling {\qococ} from C/C++., float, label=lst:qococ-usage]
    #include "qoco_custom.h"

    int main() {
       // Instantiate workspace.
       Workspace work;

       // Set default settings, but settings can be modified with 
       // work.settings.setting_name = ...
       set_default_settings(&work); 

       // Loads the default data for the problem 
       // (i.e. the data passed to the qocogen.generate() call).
       load_data(&work); 

       // Solve original problem.
       qoco_custom_solve(&work);
       printf("\nobj: %f", work.sol.obj);

       // Can modify non-zero elements of P,A,G and c, b, h.
       update_P(&work, Pnew);
       update_A(&work, Anew);
       update_G(&work, Gnew);
       update_c(&work, cnew);
       update_b(&work, bnew);
       update_h(&work, hnew);

       // Solve updated problem.
       qoco_custom_solve(&work);
       printf("\nobj: %f", work.sol.obj); 
    }
\end{lstlisting}
\FloatBarrier
\section{Numerical results}\label{sec:numerical}

We benchmark {\qoco} and {\qococ} against several conic interior-point solvers \footnote{Our benchmarks are publicly available at \url{https://github.com/qoco-org/qoco-benchmarks}}: the open-source solvers {\ecos} \cite{domahidi2013ecos} and {\clarabel} \cite{goulart2024clarabel}, the commercial solvers {\gurobi} \cite{gurobi} and {\mosek} \cite{mosek}, and against the academically licensed custom solver generator {\cvxgen} which generates a primal-dual interior point method for quadratic programs. We do not test against the custom solver generator $\textsc{bsocp}$ since it is not academically licensed. For all solvers, we use their default settings but set the tolerances $\epsilon_{abs} = \epsilon_{rel} = 10^{-7}$.

All numerical results were generated on a desktop computer with an AMD Ryzen 9 7950X3D processor running at 4.2 GHz and with 128 GB of RAM. Additionally, the model-predictive control problems in Section \ref{subsec:mpc} were also solved on the Raspberry Pi Compute Module 4 (CM4), a resource-constrained, embedded Linux system with 4GB of RAM. Our benchmarks are written in Python and we use {\cvxpy} \cite{diamond2016cvxpy,agrawal2018rewriting} to interface with the solvers. 

In our experiments, we define problem size as the total number of nonzero elements in $A$, $G$, and in the upper half of $P$. We evaluate our solvers on a set of benchmark problems we developed, a set of model-predictive control problems from \cite{kouzoupis2015towards}, the Maros–Mészáros problems \cite{maros1999repository}, and least-squares problems derived from the SuiteSparse repository \cite{davis2011university}. To evaluate the solvers' performance, we use performance profiles \cite{dolan2002benchmarking} and the shifted geometric mean, both of which are frequently used to assess solver performance \cite{goulart2024clarabel,domahidi2013ecos,osqp,garstka2021cosmo, schwan2023piqp}.

When reporting solve times for {\qococ}, we exclude code generation and compilation times and report them separately. We do this because the primary use cases for {\qocog} and {\qococ} prioritize minimizing solve time, making solve time alone the most relevant performance metric, and the one-time cost of generating and compiling the solver is less critical. Typical examples include MPC and real-time trajectory optimization, where code generation and compilation can be done offline before deploying the solvers online. 

When minimizing setup time (code generation and compilation) is important, such as iterating problem formulations or one-off solves, {\qocog} and {\qococ} may be less attractive than general-purpose solvers like {\qoco}. Additionally, as problem sizes become larger, code generation and compile times become longer and binary sizes become larger due to the explicit coding style used by custom linear algebra and depicted in Listing \ref{lst:custom-matvec}. Nevertheless, the combined code generation and compilation times remain reasonable (under ten minutes) for problems with sizes under 10,000. 

We also expect {\qococ} performance to degrade once binaries become large compared to the cache size of the processor, due to increased instruction cache misses. On desktop computers, which have larger caches, excessive code generation time and compile times are typically experienced prior to degraded performance due to instruction cache misses. On systems with comparatively smaller caches, such as the Raspberry Pi CM4, instruction cache misses start to hinder performance for problem sizes over approximately 4,000 in our benchmarks, before code generation and compilation time becomes significant.

We provide detailed numeric results in Appendix \ref{appendix:detailed-results}, which includes iteration counts and solve times for each solver on all problems we ran. Because Gurobi does not report iteration counts via the {\cvxpy} interface for SOCPs, these values are omitted from our results.

\paragraph{Shifted geometric mean}
We use the normalized shifted geometric mean to assign a scalar value to the performance of a solver on a test set of problems. The shifted geometric mean of $N$ runtimes for solver $s$ is computed as

\begin{equation*}
    g_s = \left[\prod_{p=1}^{N} (t_{s,p} + k)\right]^{1/N} - k,
\end{equation*}

where $t_{s,p}$ is the runtime of solver $s$ on problem $p$ and $k=1$ is the shift. If the solver does not find an optimal solution, then the runtime is set to 10 seconds for benchmark problems, 0.5 seconds for the MPC problems, 1200 seconds for the Maros–Mészáros problems, and 600 seconds for the SuiteSparse least-squares problems. These times were chosen to be roughly an order of magnitude higher than the slowest successful solve on the respective problem set.

We then define the normalized shifted geometric mean for solver $s$ as

\begin{equation*}
    r_s = g_s / \min_{s} g_s,
\end{equation*}

so the solver with the lowest shifted geometric mean will have a normalized shifted geometric mean of $1.00$.

\paragraph{Performance profiles}
We also plot the relative and absolute performance profiles to compare solvers. We first define the relative performance ratio of solver $s$ on problem $p$ as

\begin{equation*}
    u_{s,p} = t_{s,p} / \min_{s} t_{s, p}.
\end{equation*}

The relative performance profile then plots $f^r_s(\tau)$ where

\begin{equation*}
    f^r_s(\tau) = \frac{1}{N}\sum_{p=1}^N \mc{I}_{\leq \tau}(u_{s,p}),
\end{equation*}

and $\mc{I}_{\leq \tau}(z) = 1$ if $z \leq \tau$ and $\mc{I}_{\leq \tau}(z) = 0$ otherwise. If the relative performance profile passes through the point $(x,y)$ for solver $s$, this means that solver $s$ solves a fraction $y$ of the problems within a factor of $x$ of the fastest solver.

The absolute performance profile plots $f^a_s(\tau)$ where

\begin{equation*}
    f^a_s(\tau) = \frac{1}{N}\sum_{p=1}^N \mc{I}_{\leq \tau}(t_{s,p}).
\end{equation*}
If the absolute performance profile passes through the point $(x,y)$ for solver $s$, this means that solver $s$ solves a fraction $y$ of the problems within $x$ time. For both relative and absolute profiles, the profile of the highest-performing solver will lie above the profiles of the other solvers.

\subsection{Benchmark problems}
We solve optimization problems from five problem classes: robust Kalman filtering, group lasso regression \cite{yuan2006grouplasso}, the losslessly convexified powered-descent guidance problem \cite{acikmese2007lcvx}, Markowitz portfolio optimization \cite{markowitz1952-ds}, and the oscillating masses control problem \cite{wang2010fast}. The first three of these are SOCPs and the last two are QPs. For each of the five classes, we consider 10 different problem sizes, generating 20 unique problem instances per size, for a total of 1000 distinct optimization problems. Details on the problems we solve can be found in Appendix \ref{appendix:problem-classes}. For these problems, we plot runtime as a function of problem size in addition to reporting performance profiles and shifted geometric means. We test {\qoco} on all problems, but only test {\qococ} on the smaller problems, as code generation and compile times become prohibitively long for the larger instances.

For each solver, we solve each optimization problem 100 times and record the minimum runtime, which includes both setup and solve time. Since the computer runs various background processes that can artificially inflate execution time, we conduct 100 runs and report the minimum execution time, thereby providing a more accurate measure of the solver's performance.

Although the oscillating masses and portfolio optimization problems are QPs and can be solved with {\cvxgen}, we limit testing to the two smaller oscillating masses problems. For larger instances, {\cvxgen} was unable to generate code due to the size of the problems. For the portfolio optimization problems, we do not test {\cvxgen} due to the sparsity pattern of the factor matrix $F$ in Problem \ref{eq:portfolio}. To specify the sparsity pattern of $F$ in {\cvxgen}, we must manually provide the locations of nonzero elements into the web interface, which would have been prohibitively expensive, since $F$ has hundreds of nonzero elements. In contrast, to generate code with ${\qocog}$, we pass in problem data as sparse \texttt{scipy} matrices rather than manually specifying sparsity patterns. Although we can specify $F$ as a dense matrix and only populate the nonzero elements for {\cvxgen}, this would severely hinder its performance and result in an unfair comparison.

\begin{figure}[H]
    \captionsetup{labelfont=bf}
    \centering
    \begin{subfigure}[b]{0.49\textwidth}
        \centering
        {\includegraphics[width=\textwidth]{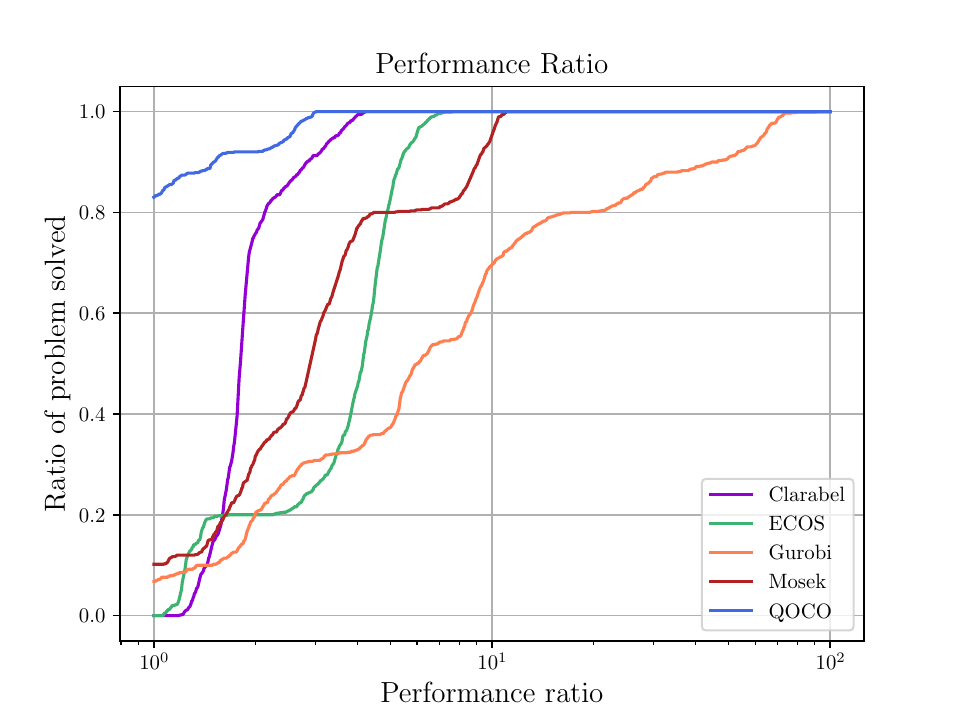}}
        \caption{Relative performance profile}
    \end{subfigure}
    \hfill
    \begin{subfigure}[b]{0.49\textwidth}
        \centering
        {\includegraphics[width=\textwidth]{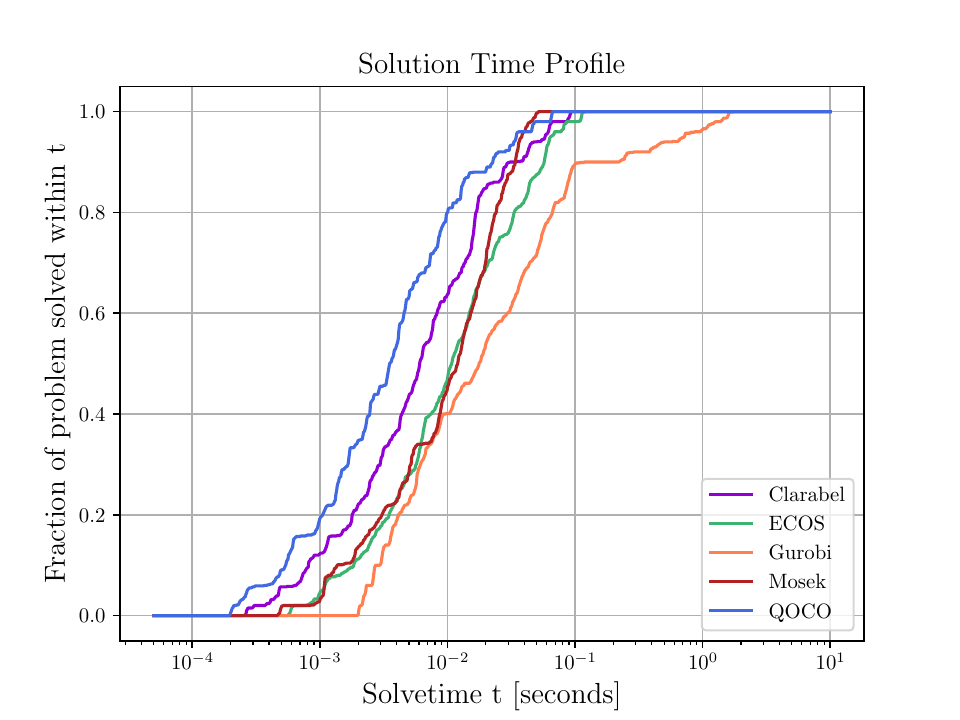}}
        \caption{Absolute performance profile}
    \end{subfigure}

    \vspace{0.5cm}

    \begin{subfigure}{1\textwidth}
        \centering
        \footnotesize
        \begin{tabular}{lccccc}
  \hline
   & \textbf{QOCO} & \textbf{Clarabel} & \textbf{ECOS} & \textbf{Gurobi} & \textbf{Mosek} \\ \hline
  Shifted GM & 1.0 & 1.6 & 2.5 & 10.2 & 1.8 \\ 
  Failure Rate (\%) & 0.0 & 0.0 & 0.0 & 0.0 & 0.0 \\ \hline 
\end{tabular}

        \caption{Shifted geometric means and failure rates}
      \end{subfigure}
    \caption{\bf Performance profiles for all benchmark problems}
    \label{fig:bench}
\end{figure}

\begin{figure}[H]
    \captionsetup{labelfont=bf}
    \centering
    \begin{subfigure}[b]{0.49\textwidth}
        \centering
        {\includegraphics[width=\textwidth]{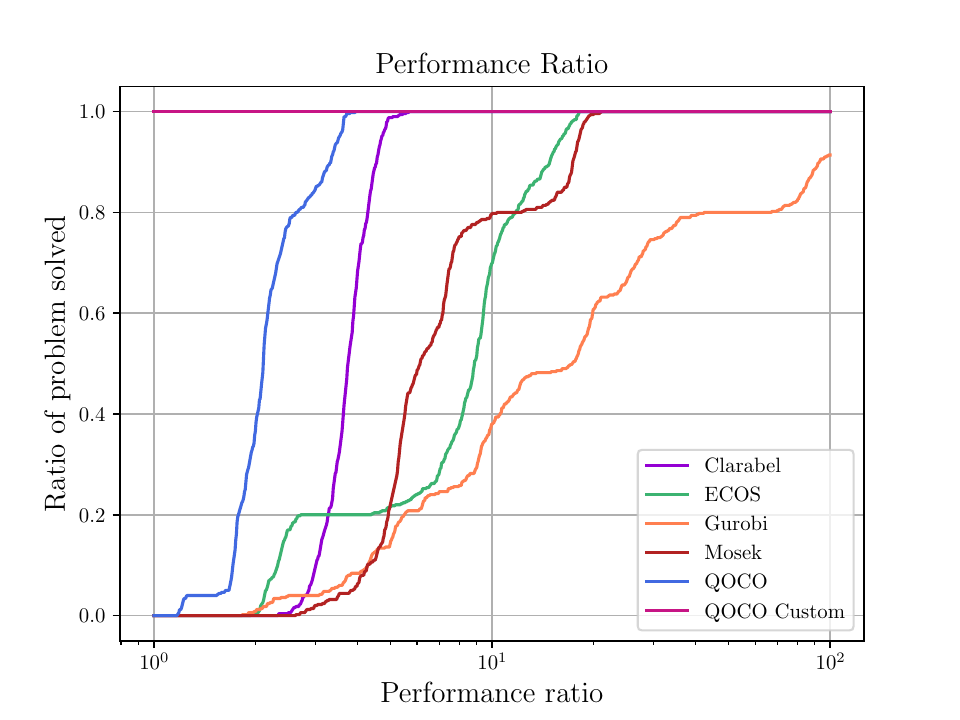}}
        \caption{Relative performance profile}
    \end{subfigure}
    \hfill
    \begin{subfigure}[b]{0.49\textwidth}
        \centering
        {\includegraphics[width=\textwidth]{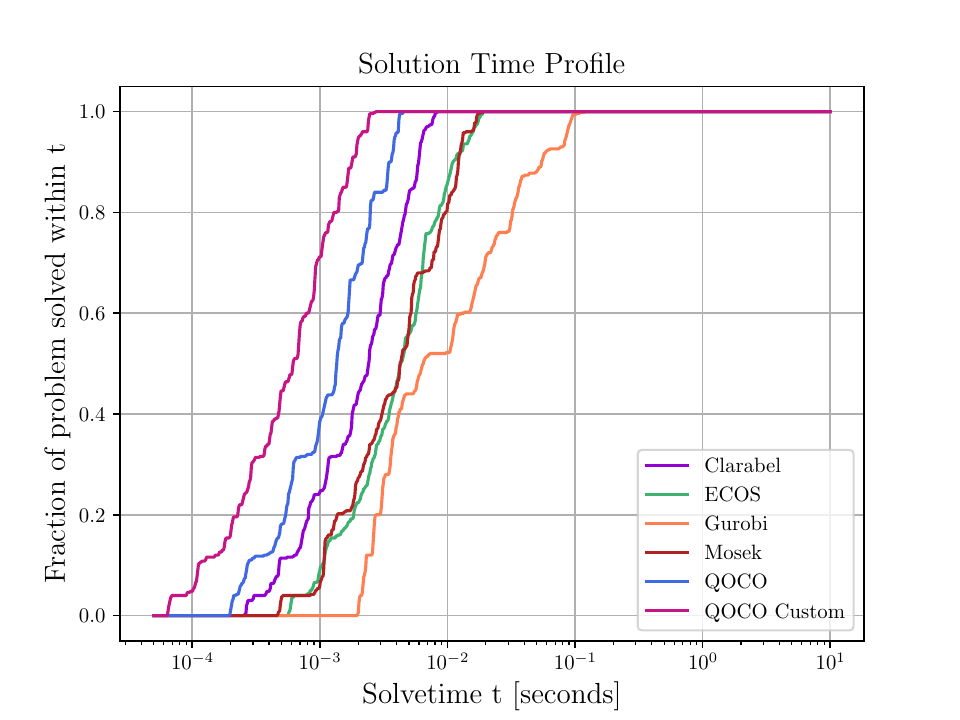}}
        \caption{Absolute performance profile}
    \end{subfigure}

    \vspace{0.5cm}

    \begin{subfigure}{1\textwidth}
        \centering
        \footnotesize
        \begin{tabular}{lcccccc}
  \hline
    & \textbf{QOCO Custom}   & \textbf{QOCO} & \textbf{Clarabel} & \textbf{ECOS} & \textbf{Gurobi} & \textbf{Mosek} \\ \hline
  Shifted GM & 1.0 & 2.0 & 3.5 & 6.9 & 23.9 & 7.0 \\ 
  Failure Rate (\%) & 0.0 & 0.0 & 0.0 & 0.0 & 0.0 & 0.0 \\ \hline 
\end{tabular}

        \caption{Shifted geometric means and failure rates}
      \end{subfigure}
    \caption{\bf Performance profiles for benchmark problems with custom solver}
    \label{fig:bench-custom}
\end{figure}

When computing normalized shifted geometric mean, relative performance profiles, and absolute performance profiles, all solvers must be tested on the same set of problems. Since we test $\qococ$ on only half of the problems we provide two versions of the aforementioned metrics. The first version, given by Figure \ref{fig:bench}, includes {\clarabel}, {\ecos}, {\gurobi}, {\mosek}, and {\qoco} evaluated on all problems. The second version, given by Figure \ref{fig:bench-custom}, also includes {\qococ}, but the solvers are evaluated on the half of the benchmark problems which {\qococ} is tested on.

\begin{figure}[H]
    \captionsetup{labelfont=bf}
    \centering
    \includegraphics[width=\textwidth]{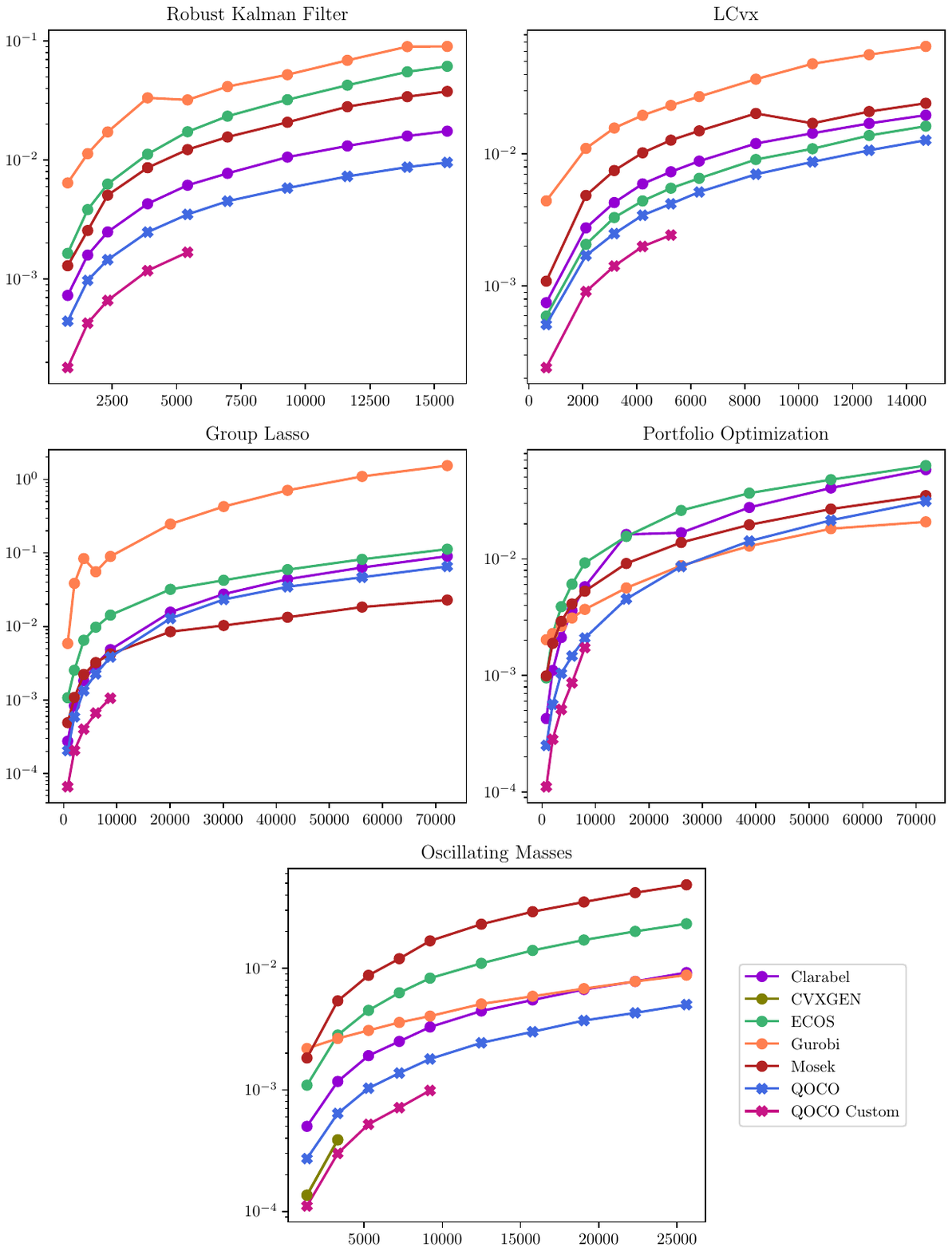}
    \caption{\bf Solvetime in seconds vs problem size for benchmark problems}
    \label{fig:bench-time-vs-size}
\end{figure}

Figures \ref{fig:bench-time-vs-size}, \ref{fig:bench}, and \ref{fig:bench-custom} show that {\qococ} is the fastest solver on our benchmark problems by a significant margin and is as performant as {\cvxgen} on the $40$ smallest oscillating mass problems. This latter result is unsurprising, as {\qococ} uses an identical algorithm to {\cvxgen} when solving QPs and both utilize custom linear algebra routines. Additionally, {\qoco} is the fastest open-source generic solver on all of the problems, although it is slower than Mosek on the group lasso problems and slower than Gurobi on some of the largest portfolio optimization problems. However, we notice that if we disable presolve for Mosek and Gurobi, {\qoco} is faster on these problems as well.

Table \ref{tab:codegen_stats_benchmark} reports the code generation time, compilation time, code size, and binary size of the generated solvers. We see that all problems have sizes under 10,000 and as the problem size increases, code generation time, compilation time, code size, and binary size all increase, but all solvers still take less than eight minutes to generate and compile.

\begin{table}[ht]
\centering
\footnotesize
\begin{tabular}{lrrrrr}
\toprule
Problem & Size & Codegen Time (s) & Compile Time (s) & Code Size (KB) & Binary Size (KB) \\
\midrule
group\_lasso\_N\_1 & 761 & 4.2 & 4.3 & 545 & 234 \\
group\_lasso\_N\_2 & 2022 & 16.9 & 77.0 & 1292 & 590 \\
group\_lasso\_N\_3 & 3783 & 38.7 & 37.6 & 2397 & 1150 \\
group\_lasso\_N\_4 & 6044 & 70.5 & 69.3 & 3923 & 1898 \\
group\_lasso\_N\_5 & 8805 & 112.3 & 96.9 & 5949 & 2918 \\
\midrule
portfolio\_N\_2 & 804 & 2.8 & 4.1 & 481 & 226 \\
portfolio\_N\_4 & 2008 & 11.1 & 45.0 & 1082 & 518 \\
portfolio\_N\_6 & 3612 & 25.3 & 473.9 & 1890 & 918 \\
portfolio\_N\_8 & 5616 & 45.5 & 46.0 & 2982 & 1482 \\
portfolio\_N\_10 & 8020 & 72.9 & 106.1 & 6309 & 3026 \\
\midrule
lcvx\_N\_15 & 644 & 3.6 & 19.3 & 1057 & 482 \\
lcvx\_N\_50 & 2114 & 39.3 & 57.0 & 3795 & 1590 \\
lcvx\_N\_75 & 3164 & 89.9 & 85.5 & 5694 & 2434 \\
lcvx\_N\_100 & 4214 & 163.8 & 140.1 & 7934 & 3322 \\
lcvx\_N\_125 & 5264 & 263.0 & 193.3 & 9675 & 4026 \\
\midrule
robust\_kalman\_filter\_N\_25 & 775 & 5.8 & 11.3 & 759 & 318 \\
robust\_kalman\_filter\_N\_50 & 1550 & 23.2 & 128.8 & 1491 & 582 \\
robust\_kalman\_filter\_N\_75 & 2325 & 53.0 & 34.1 & 2236 & 890 \\
robust\_kalman\_filter\_N\_125 & 3875 & 155.2 & 76.8 & 3730 & 1450 \\
robust\_kalman\_filter\_N\_175 & 5425 & 311.0 & 122.5 & 5235 & 2022 \\
\midrule
oscillating\_masses\_N\_8 & 1344 & 4.7 & 22.9 & 1112 & 478 \\
oscillating\_masses\_N\_20 & 3312 & 28.2 & 33.3 & 2785 & 1150 \\
oscillating\_masses\_N\_32 & 5280 & 72.4 & 60.4 & 4450 & 1806 \\
oscillating\_masses\_N\_44 & 7248 & 138.5 & 86.3 & 6122 & 2450 \\
oscillating\_masses\_N\_56 & 9216 & 226.9 & 123.6 & 7839 & 3102 \\
\bottomrule
\end{tabular}
\captionsetup{labelfont=bf}
\caption{ \bf Code generation time, compilation time, and resulting code and binary sizes for benchmark problems.}
\label{tab:codegen_stats_benchmark}
\end{table}

\subsection{Model-predictive control problems}\label{subsec:mpc}
We consider a set of 64 model-predictive control problems taken from \cite{kouzoupis2015towards} and also used by \cite{goulart2024clarabel}. Note that we exclude the eight ``nonlinear chain'' problems as they took too long to generate custom solvers for. We run these problems on our desktop computer with the AMD Ryzen 9 7950X3D, as well as the Raspberry Pi CM4. Similarly to the benchmark problems, we solve each MPC problem 100 times and record the minimum runtime. These problems take the following form

\begin{equation*}
    \begin{split}
        \underset{x, y, u}{\mathrm{minimize}} 
        \quad & \sum_{k = 0}^{\tau} \begin{bmatrix}y_k - y_k^r \\ u_k - u_k^r\end{bmatrix}^\top 
        \begin{bmatrix}Q & S \\ S^\top & R \end{bmatrix}
        \begin{bmatrix}y_k - y_k^r \\ u_k - u_k^r\end{bmatrix} + (x_\tau - x_\tau^r)^\top P (x_\tau - x_\tau^r) \\
        \mathrm{subject\;to} 
        \quad & x_0 = x_{\mathrm{init}} \\
        & x_{k+1} = Ax_k + Bu_k + f_k \\
        & y_k = Cx_k + Du_k + e_k \\
        & d_{\min} \le Mx_k + Nu_k \le d_{\max} \\
        & u_{\min} \le u_k \le u_{\max} \\
        & y_{\min} \le y_k \le y_{\max} \\
        & d^\tau_{\min} \le Tx_\tau \le d^\tau_{\max}. \\
    \end{split}
\end{equation*}

Table \ref{tab:codegen_stats_mpc} shows that all problems except \texttt{SpringMass\_1} problem have sizes under 10,000 and combined generation and compilation times under ten minutes. Additionally, among problems with sizes under 10,000, all but \texttt{quadcopter\_3} require less than five minutes to generate and compile solvers.

\begin{figure}[H]
    \captionsetup{labelfont=bf}
    \centering
    \begin{subfigure}[b]{0.49\textwidth}
        \centering
        {\includegraphics[width=\textwidth]{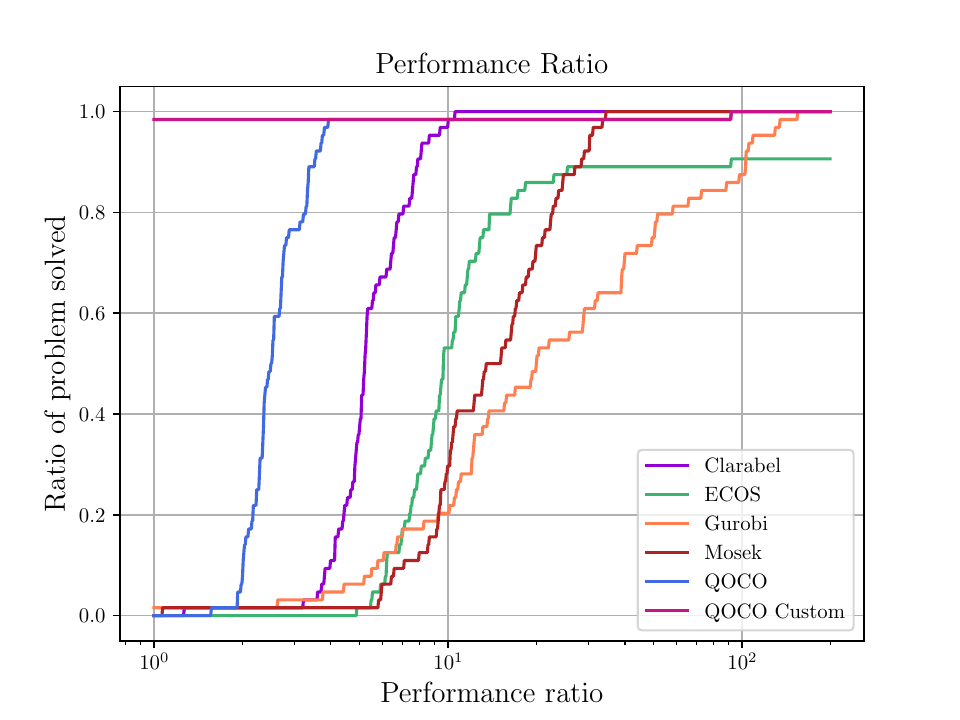}}
        \caption{Relative performance profile}
    \end{subfigure}
    \hfill
    \begin{subfigure}[b]{0.49\textwidth}
        \centering
        {\includegraphics[width=\textwidth]{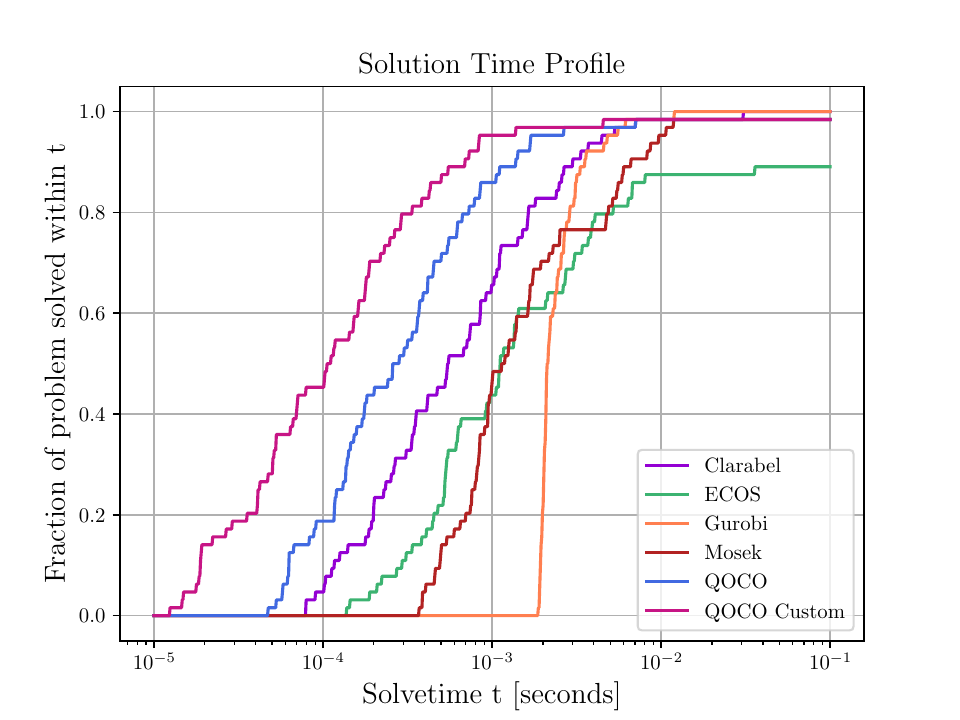}}
        \caption{Absolute performance profile}
    \end{subfigure}

    \vspace{0.5cm}

    \begin{subfigure}{1\textwidth}
        \centering
        \footnotesize
        \begin{tabular}{lcccccc}
  \hline
    & \textbf{QOCO Custom}   & \textbf{QOCO} & \textbf{Clarabel} & \textbf{ECOS} & \textbf{Gurobi} & \textbf{Mosek} \\ \hline
  Shifted GM & 4.4 & 4.5 & 1.0 & 31.4 & 1.7 & 2.5 \\ 
  Failure Rate (\%) & 1.6 & 1.6 & 0.0 & 10.9 & 0.0 & 0.0 \\ \hline 
\end{tabular}

        \caption{Shifted geometric means and failure rates}
      \end{subfigure}
    \caption{\bf Performance profiles for model-predictive control problems}
    \label{fig:mpc}
\end{figure}

Figure \ref{fig:mpc} shows the numerical results on the desktop computer. We see that {\clarabel}, {\gurobi}, and {\mosek} solve all problems successfully, {\qoco} and {\qococ} each fail on only one, and {\ecos} fails on seven. Among the 63 problems solved by both {\qoco} and {\qococ}, {\qococ} is the fastest solver by a wide margin, followed by {\qoco}.

\begin{figure}[H]
    \captionsetup{labelfont=bf}
    \centering
    \begin{subfigure}[b]{0.49\textwidth}
        \centering
        {\includegraphics[width=\textwidth]{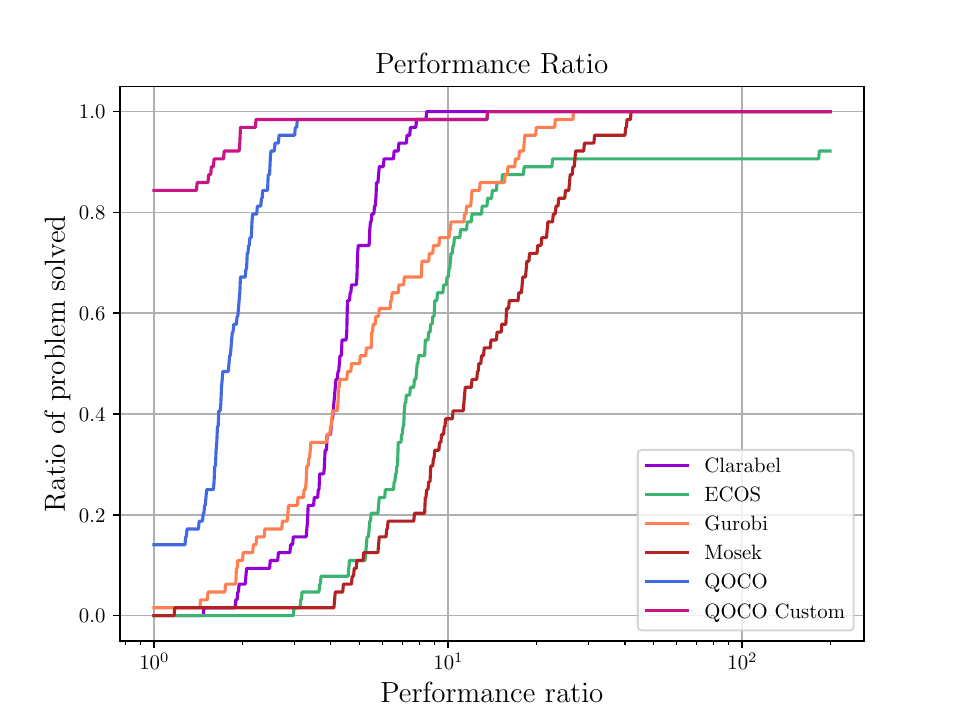}}
        \caption{Relative performance profile}
    \end{subfigure}
    \hfill
    \begin{subfigure}[b]{0.49\textwidth}
        \centering
        {\includegraphics[width=\textwidth]{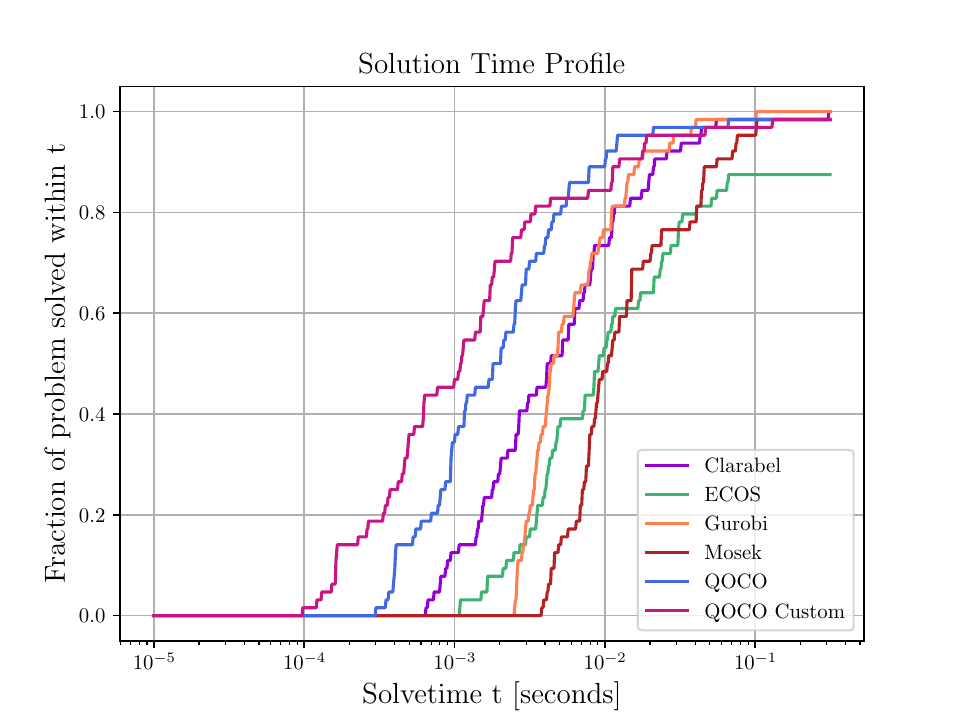}}
        \caption{Absolute performance profile}
    \end{subfigure}

    \vspace{0.5cm}

    \begin{subfigure}{1\textwidth}
        \centering
        \footnotesize
        \begin{tabular}{lcccccc}
  \hline
    & \textbf{QOCO Custom}   & \textbf{QOCO} & \textbf{Clarabel} & \textbf{ECOS} & \textbf{Gurobi} & \textbf{Mosek} \\ \hline
  Shifted GM & 1.3 & 1.1 & 1.4 & 7.0 & 1.0 & 3.6 \\ 
  Failure Rate (\%) & 1.6 & 1.6 & 0.0 & 10.9 & 0.0 & 0.0 \\ \hline 
\end{tabular}

        \caption{Shifted geometric means and failure rates}
      \end{subfigure}
    \caption{\bf Performance profiles for model-predictive control problems on Raspberry Pi CM4}
    \label{fig:mpc_cm4}
\end{figure}

Figure \ref{fig:mpc_cm4} shows the results on the Raspberry Pi CM4, and we see that all solvers solve the same number of problems as on the desktop. We also see that problems are solved roughly an order magnitude slower on the CM4, as it is a less powerful processor. We again see that among the 63 problems solved by both {\qoco} and {\qococ}, {\qococ} or {\qoco} is always the fastest solver, but on the nine largest problems, {\qoco} is faster than {\qococ}. This occurs for the largest problems, since the {\qococ} binaries are larger than the CM4's L2 cache, likely leading to many instruction cache misses and increased solve time.  

We observe {\qoco} outperforming {\qococ} only on the CM4 and not on the desktop computer. This is due to the substantially smaller cache on the CM4 (1 MB of L2 cache) compared to the AMD Ryzen 9 7950X3D (16 MB of L2 cache and 128 MB of L3 cache). On the desktop processor, the {\qococ} binaries for all tested problems fit comfortably within the L2 cache. If we were to benchmark problems on the desktop where the {\qococ} binary approached the 16 MB L2 cache size, we would expect to see some performance degradation. If the binary size approached the 128 MB L3 cache limit, we expect the performance degradation to be significant and could lead to {\qoco} outperforming {\qococ}. However, for problems of this scale, code generation and compilation time will be prohibitively long, on the order of several hours to even a day, making the use of {\qocog} and {\qococ} impractical. 

{\footnotesize
\begin{longtable}{lrrrrr}
\caption{\bf Code generation time, compilation time, and resulting code and binary sizes for mpc problem instances.}
\label{tab:codegen_stats_mpc} \\

\toprule
Problem & Size & Codegen Time (s) & Compile Time (s) & Code Size (KB) & Binary Size (KB) \\
\midrule
\endfirsthead

\toprule
Problem & Size & Codegen Time (s) & Compile Time (s) & Code Size (KB) & Binary Size (KB) \\
\midrule
\endhead

\midrule
\multicolumn{6}{r}{\footnotesize Continued on next page} \\
\endfoot

\bottomrule
\endlastfoot
aircraft\_1 & 504 & 1.4 & 2.7 & 412 & 182 \\
aircraft\_2 & 524 & 1.4 & 2.6 & 413 & 186 \\
aircraft\_3 & 464 & 1.1 & 2.3 & 376 & 170 \\
aircraft\_4 & 504 & 1.4 & 2.7 & 412 & 182 \\
aircraft\_10 & 584 & 2.2 & 3.8 & 470 & 198 \\
aircraft\_11 & 636 & 2.3 & 3.9 & 480 & 202 \\
aircraft\_12 & 636 & 2.3 & 3.8 & 480 & 202 \\
aircraft\_13 & 3116 & 56.0 & 32.7 & 2230 & 890 \\
\midrule
ballOnPlate\_1 & 398 & 1.5 & 2.4 & 328 & 150 \\
ballOnPlate\_2 & 398 & 1.5 & 2.4 & 328 & 150 \\
ballOnPlate\_3 & 398 & 1.5 & 2.4 & 328 & 150 \\
ballOnPlate\_4 & 658 & 4.1 & 4.3 & 509 & 214 \\
\midrule
binaryDistillationColumn\_1 & 2936 & 8.9 & 179.2 & 3107 & 1338 \\
binaryDistillationColumn\_2 & 2936 & 8.8 & 179.4 & 3107 & 1338 \\
\midrule
dcMotor\_1 & 564 & 1.8 & 3.4 & 456 & 194 \\
dcMotor\_2 & 1114 & 6.7 & 18.1 & 858 & 354 \\
dcMotor\_3 & 5514 & 175.3 & 78.4 & 4184 & 1602 \\
dcMotor\_4 & 5514 & 175.3 & 78.5 & 4184 & 1602 \\
dcMotor\_5 & 1114 & 6.7 & 18.0 & 858 & 354 \\
dcMotor\_6 & 1114 & 6.8 & 17.8 & 857 & 354 \\
\midrule
doubleInvertedPendulum\_1 & 502 & 1.1 & 2.4 & 390 & 178 \\
doubleInvertedPendulum\_2 & 502 & 1.2 & 2.4 & 389 & 178 \\
doubleInvertedPendulum\_3 & 542 & 1.4 & 2.7 & 425 & 190 \\
\midrule
fiordosExample\_1 & 119 & 0.1 & 0.6 & 132 & 74 \\
fiordosExample\_2 & 139 & 0.2 & 0.7 & 150 & 78 \\
fiordosExample\_3 & 119 & 0.1 & 0.6 & 132 & 74 \\
\midrule
forcesExample\_1 & 171 & 0.2 & 0.8 & 163 & 86 \\
forcesExample\_2 & 170 & 0.2 & 0.8 & 161 & 86 \\
forcesExample\_3 & 189 & 0.3 & 0.8 & 174 & 90 \\
forcesExample\_4 & 369 & 1.1 & 1.8 & 287 & 134 \\
\midrule
helicopter\_1 & 201 & 0.3 & 0.9 & 197 & 98 \\
helicopter\_2 & 1051 & 7.7 & 14.0 & 863 & 378 \\
helicopter\_3 & 536 & 2.1 & 3.2 & 449 & 186 \\
\midrule
nonlinearCstr\_1 & 1164 & 7.2 & 16.3 & 829 & 338 \\
nonlinearCstr\_2 & 1164 & 7.2 & 16.5 & 829 & 338 \\
nonlinearCstr\_3 & 294 & 0.5 & 1.3 & 247 & 114 \\
\midrule
pendulum\_1 & 435 & 1.1 & 2.0 & 332 & 154 \\
pendulum\_2 & 432 & 1.1 & 2.0 & 329 & 150 \\
pendulum\_3 & 155 & 0.1 & 0.7 & 153 & 86 \\
\midrule
quadcopter\_1 & 1592 & 15.0 & 118.2 & 1770 & 782 \\
quadcopter\_2 & 3172 & 60.0 & 57.0 & 3656 & 1594 \\
quadcopter\_3 & 7912 & 409.8 & 197.9 & 9408 & 3974 \\
quadcopter\_4 & 3172 & 60.0 & 57.0 & 3655 & 1594 \\
quadcopter\_5 & 3072 & 51.3 & 55.1 & 3686 & 1590 \\
quadcopter\_6 & 1678 & 15.7 & 123.5 & 1797 & 794 \\
\midrule
robotArm\_1 & 276 & 0.6 & 1.3 & 243 & 114 \\
robotArm\_2 & 1056 & 8.9 & 12.7 & 795 & 330 \\
\midrule
shell\_1 & 739 & 4.1 & 15.7 & 974 & 458 \\
shell\_2 & 1469 & 16.0 & 146.7 & 2017 & 858 \\
shell\_3 & 739 & 4.0 & 15.7 & 974 & 458 \\
\midrule
spacecraft\_1 & 1137 & 8.6 & 15.6 & 835 & 338 \\
spacecraft\_2 & 1137 & 8.6 & 15.8 & 835 & 338 \\
\midrule
springMass\_1 & 21239 & 2521.6 & 567.3 & 17247 & 6326 \\
springMass\_2 & 2159 & 20.4 & 178.8 & 1704 & 658 \\
springMass\_3 & 2144 & 20.3 & 167.4 & 1690 & 654 \\
springMass\_4 & 4279 & 82.3 & 53.5 & 3382 & 1314 \\
\midrule
toyExample\_1 & 199 & 0.3 & 0.8 & 178 & 90 \\
toyExample\_2 & 389 & 1.1 & 1.9 & 296 & 138 \\
toyExample\_3 & 199 & 0.3 & 0.8 & 178 & 90 \\
toyExample\_4 & 959 & 6.5 & 8.2 & 650 & 278 \\
toyExample\_5 & 1909 & 26.4 & 99.8 & 1259 & 490 \\
\midrule
tripleInvertedPendulum\_1 & 2739 & 22.2 & 272.7 & 2109 & 822 \\
tripleInvertedPendulum\_2 & 2919 & 28.2 & 28.3 & 2291 & 926 \\
tripleInvertedPendulum\_3 & 2919 & 28.1 & 28.3 & 2289 & 926 \\
\end{longtable}
}

\subsection{Maros–Mészáros problems}
The Maros–Mészáros problems are a set of 138 challenging QPs that include a wide range of problem sizes (from $3$ to around $300,000$ variables) and contain very difficult problems that have ill-conditioning and poor scaling \cite{maros1999repository}. These properties make the Maros–Mészáros problems a good test set to assess the robustness of solvers. We test {\qoco} but not {\qococ} on these problems since many problems are far too large to generate code for. Additionally, we run each problem once rather than the 100 runs we do for the benchmark and MPC problems since many of the Maros–Mészáros problems are extremely large, and running them 100 times per solver is extremely cumbersome.

\begin{figure}[H]
    \captionsetup{labelfont=bf}
    \centering
    \begin{subfigure}[b]{0.49\textwidth}
        \centering
        {\includegraphics[width=\textwidth]{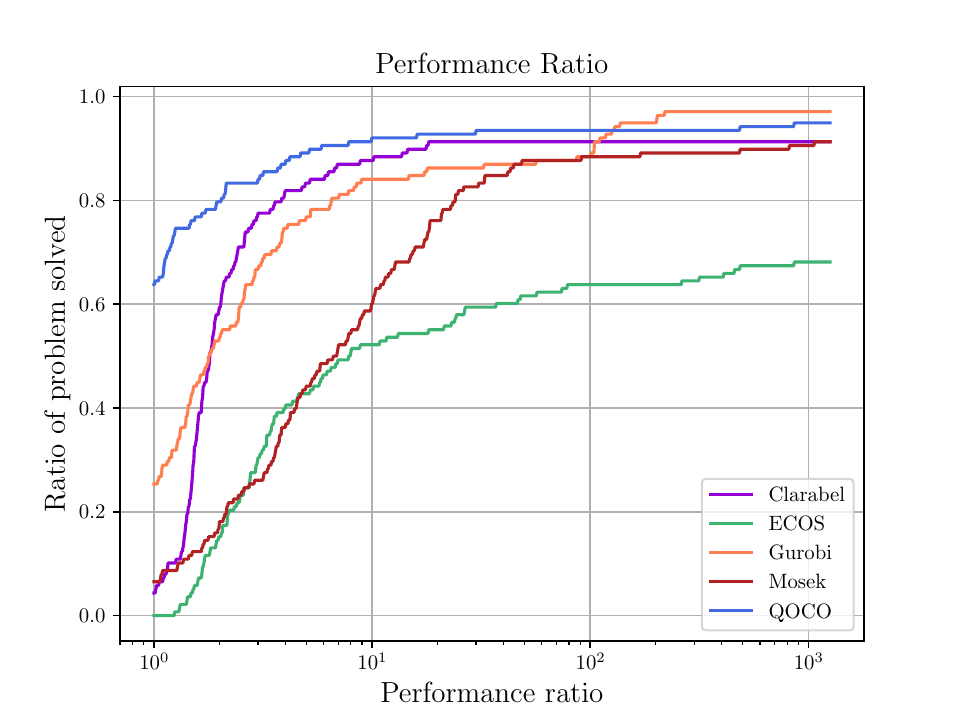}}
        \caption{Relative performance profile}
    \end{subfigure}
    \hfill
    \begin{subfigure}[b]{0.49\textwidth}
        \centering
        {\includegraphics[width=\textwidth]{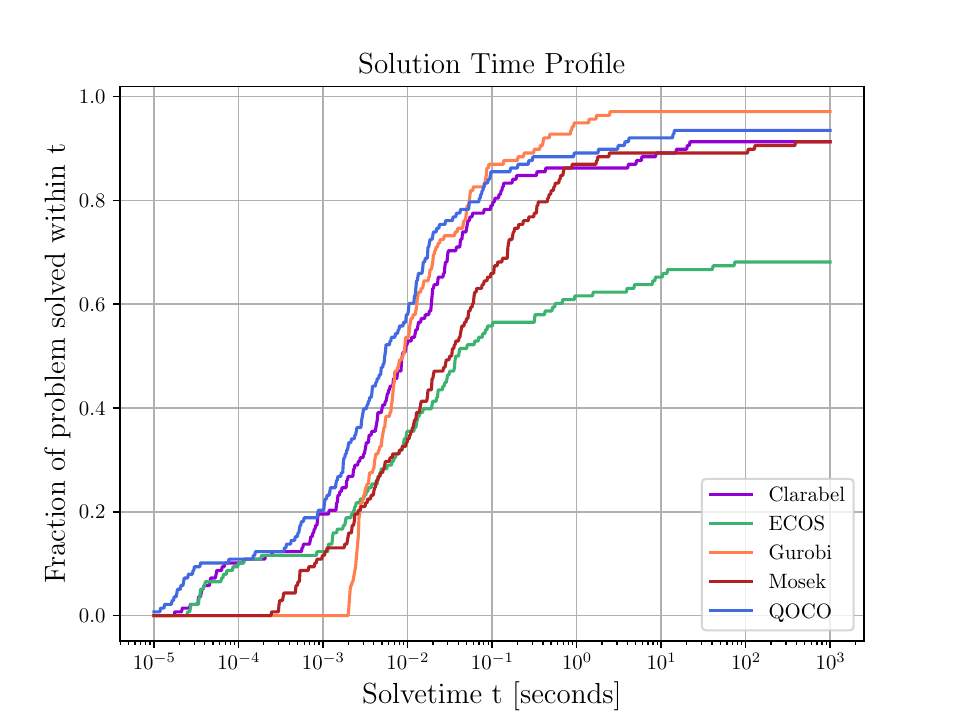}}
        \caption{Absolute performance profile}
    \end{subfigure}

    \vspace{0.5cm}

    \begin{subfigure}{1\textwidth}
        \centering
        \footnotesize
        \begin{tabular}{lccccc}
  \hline
   & \textbf{QOCO} & \textbf{Clarabel} & \textbf{ECOS} & \textbf{Gurobi} & \textbf{Mosek} \\ \hline
  Shifted GM & 2.4 & 3.6 & 33.9 & 1.0 & 4.1 \\ 
  Failure Rate (\%) & 6.5 & 8.7 & 31.9 & 2.9 & 8.7 \\ \hline 
\end{tabular}

        \caption{Shifted geometric means and failure rates}
      \end{subfigure}
    \caption{\bf Performance profiles for Maros–Mészáros problems}
    \label{fig:maros}
\end{figure}

Figure \ref{fig:maros} shows that the proprietary solver Gurobi successfully solved the largest fraction of the problems followed by {\qoco}, {\clarabel}, {\mosek}, and {\ecos}. We can also see that for the majority of the problems, {\qoco} was the fastest solver, whereas {\mosek} and {\ecos} tended to be the slowest. The latter result is unsurprising as {\mosek} and {\ecos} can only handle linear objectives and the reformulation of the quadratic objective into a second-order cone likely introduces significant fill-in, slowing these solvers down. We also observe that for many of the largest problems, Gurobi is the fastest solver, due to multithreading in its matrix factorization.

\subsection{SuiteSparse least-squares problems}
Finally, we consider a set of 23 least-squares problems $Ax \approx b$, where the matrix $A \in \R^{m \times n}$ is drawn from the SuiteSparse Matrix Collection \cite{davis2011university}. Following \cite{osqp} and \cite{goulart2024clarabel}, we get an approximate solution to each least-squares problem by solving a Huber regression and lasso regression problem for a total of 46 problems, which can both be formulated as constrained quadratic programs \cite{Boyd2004}. These optimization problems are quite large with some having over a million optimization variables, so we do not test {\qococ} due to the prohibitively long code-generation time. Additionally, we run each problem once rather than 100 times.

The Huber regression problem is given by
\begin{equation*}\label{eq:huber}
    \begin{split}
        \underset{x}{\text{minimize}} 
        \quad & \sum_{i=1}^{m} \phi(a_i^\top x - b_i),  \\
    \end{split}
\end{equation*}

where $a_i^\top$ is the $i^{th}$ row of $A$ and $\phi: \R \to \R$ is defined as

\begin{equation*}
    \phi(z) = 
    \begin{cases}
    z^2 & \text{if } |z| \leq \delta \\
    \delta (2|z| - \delta) & \text{if } |z| > \delta
    \end{cases},
\end{equation*}

where we use $\delta=1$.

The lasso regression problem is given by

\begin{equation*}\label{eq:lasso}
    \begin{split}
        \underset{x}{\text{minimize}} 
        \quad & \|Ax - b\|_2^2 + \lambda \|x\|_1,  \\
    \end{split}
\end{equation*}

where we set $\lambda = \|A^\top b\|_\infty$. 

\begin{figure}[H]
    \captionsetup{labelfont=bf}
    \centering
    \begin{subfigure}[b]{0.49\textwidth}
        \centering
        {\includegraphics[width=\textwidth]{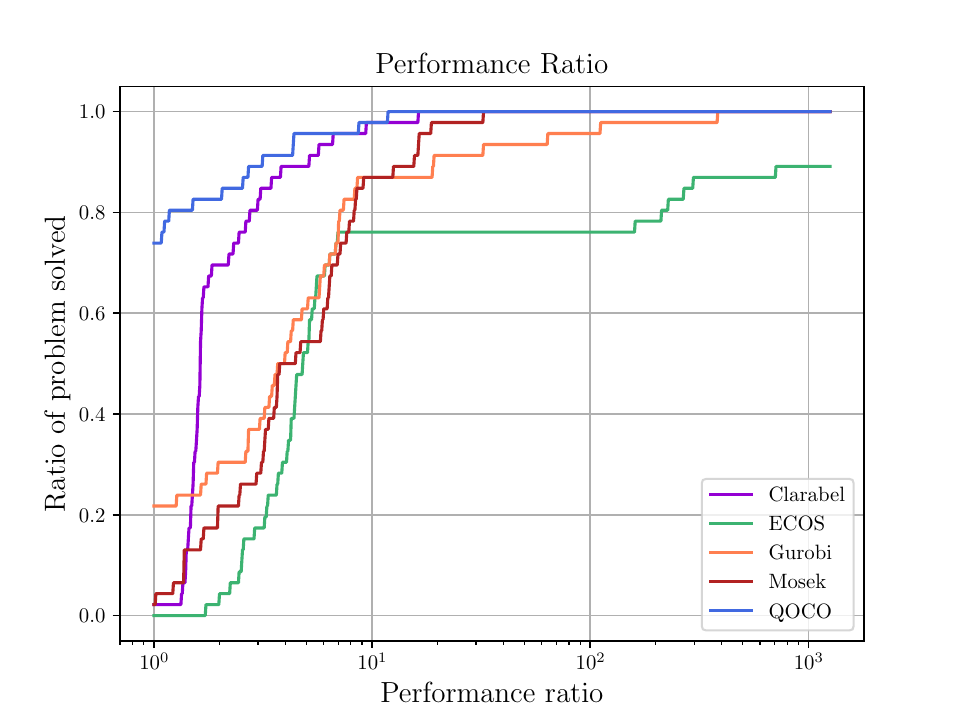}}
        \caption{Relative performance profile}
    \end{subfigure}
    \hfill
    \begin{subfigure}[b]{0.49\textwidth}
        \centering
        {\includegraphics[width=\textwidth]{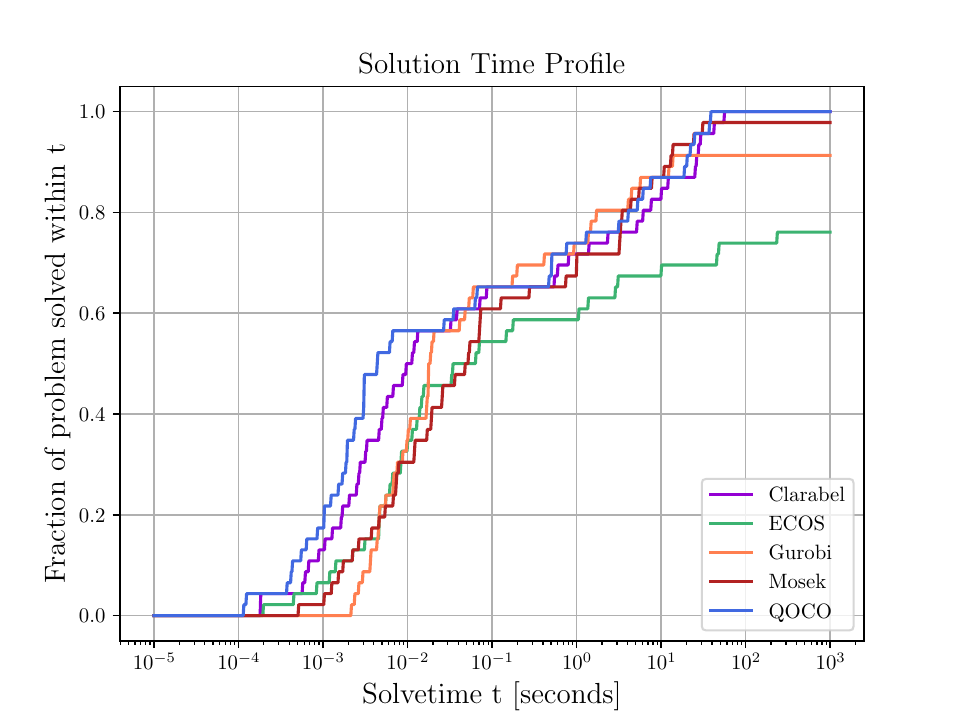}}
        \caption{Absolute performance profile}
    \end{subfigure}

    \vspace{0.5cm}

    \begin{subfigure}{1\textwidth}
        \centering
        \footnotesize
        \begin{tabular}{lccccc}
  \hline
   & \textbf{QOCO} & \textbf{Clarabel} & \textbf{ECOS} & \textbf{Gurobi} & \textbf{Mosek} \\ \hline
  Shifted GM & 1.0 & 1.2 & 7.6 & 1.6 & 1.2 \\ 
  Failure Rate (\%) & 0.0 & 0.0 & 23.9 & 8.7 & 2.2 \\ \hline 
\end{tabular}

        \caption{Shifted geometric means and failure rates}
      \end{subfigure}
    \caption{\bf Performance profiles for SuiteSparse least-squares problems}
    \label{fig:suitesparse}
\end{figure}

Figure \ref{fig:suitesparse} shows that {\qoco} and {\clarabel} successfully solve all problems, with {\qoco} being the fastest on most instances. After {\qoco} and {\clarabel}, the solvers which solve the largest fraction of problems are {\mosek}, {\gurobi}, and {\ecos} in that order. As with the Maros–Mészáros problems, we observe that Gurobi is the fastest on some of the largest problems, due to it multithreaded matrix factorization.

In total, there are 184 Maros–Mészáros problems and SuiteSparse least-squares problems, and {\gurobi} solves 176, {\qoco} solves 175, {\clarabel} solves 172, {\mosek} solves 171, and {\ecos} solves 129.

\section{Conclusion}

We have presented {\qoco}, a C-based solver, and {\qocog}, a custom solver generator for quadratic objective second-order cone programs (SOCPs). Both implement primal-dual interior-point methods, with {\qocog} generating custom solvers written in C, called {\qococ}, that exploit the sparsity pattern of problems to gain a computational advantage. We demonstrate that {\qoco} is a fast and robust solver, and {\qococ} is significantly faster than {\qoco}, making it a useful tool for real-time applications, such as model-predictive control, trajectory optimization, and other domains where computational efficiency is critical. 

Since both {\qoco} and {\qocog} are open-source, they are accessible to users in both academia and industry. {\qoco} and {\qocog} are integrated with {\cvxpy} and {\cpg} respectively, allowing users to formulate optimization problems in a natural way following from math rather than manually converting the problem to the solver-required standard form. This integration enhances usability of our solvers. 

Documentation and installation instructions for {\qoco} and {\qocog} are available at \url{https://qoco-org.github.io/qoco/index.html}.

\paragraph{Acknowledgments}
This research was supported by Blue Origin LLC and ONR grants N000142512231 and N00014-25-1-2319. Government sponsorship is acknowledged. The authors would like to thank Danylo Malyuta for his review of this paper.

\newpage
\bibliography{references}
\appendix 
\section{Nonsymmetric Newton system}\label{appendix:asymmetric-system}
In this section we will show that Newton steps applied to Equations \eqref{eq:central-path-1}, \eqref{eq:central-path-2}, \eqref{eq:central-path-3}, \eqref{eq:central-path-4} result in a nonsymmetric linear system which is less efficient to store and factor than the symmetric system presented in Equation \eqref{eq:kkt-system}.

We must first introduce the arrow matrix, a matrix that will allow us to rewrite the Jordan product, defined in Equation \eqref{eq:jordan-product}, as a matrix-vector multiplication. The arrow matrix can be written for vectors in the non-negative orthant and second-order cone as

\begin{equation*}
    \mathrm{Arw}(x) =
    \begin{cases}
        \mathrm{diag}(x) \qquad\;\; \text{if } x \in \R^l_+ \\
        \begin{bmatrix}
            x_0 & x_1^\top \\
            x_1 & x_0I
        \end{bmatrix} \quad \text{if } x \in \mc{Q}^q.
        \end{cases}    
\end{equation*}

For a vector in cone $\mc{K}$ where 

\begin{equation*}
    \mc{K} = \mc{C}_1 \times \mc{C}_2 \times \cdots \times \mc{C}_K,
\end{equation*}

and $\mc{C}_i$ is the non-negative orthant or a second-order cone, we can write the arrow matrix as

\begin{equation*}
    \mathrm{Arw}(x) = \mathrm{blkdiag}(\mathrm{Arw}(x_1), \ldots, \mathrm{Arw}(x_K)) \quad \text{where } x_i \in \mc{C}_i,
\end{equation*}

Note that if $x$ is in the interior of $\mc{K}$, then $\mathrm{Arw}(x)$ is positive definite and thus is invertible.

We can then write the Jordan product for two vectors $u,v \in \mc{K}$ as

\begin{equation*}
    x \circ y = \mathrm{Arw}(x)y,
\end{equation*}

We now linearize Equations \eqref{eq:symm-central-path-1} - \eqref{eq:symm-central-path-5} about the current iterate, $(x_{k},s_{k},y_{k},z_{k})$, write the Jordan products using the arrow matrix and get the linear system

\begin{equation*}
    \begin{bmatrix}
        P & 0 & A^\top & G^\top \\
        0 & \mathrm{Arw}(z_k) & 0 & \mathrm{Arw}(s_k) \\
        A & 0 & 0 & 0 \\
        G & I & 0 & 0
    \end{bmatrix}
    \begin{bmatrix}
        \Delta x \\
        \Delta s \\
        \Delta y \\
        \Delta z
    \end{bmatrix}
    =
    \begin{bmatrix}
        -r_x \\
        -r_s \\
        -r_y \\
        -r_z
    \end{bmatrix}
\end{equation*}

We can attempt to symmetrize the coefficient matrix by multiplying the second equation by $\mathrm{Arw}(s_k)^{-1}$ to get

\begin{equation}\label{eq:asymmetric-system}
    \begin{bmatrix}
        P & 0 & A^\top & G^\top \\
        0 & \mathrm{Arw}(s_k)^{-1}\mathrm{Arw}(z_k) & 0 & I \\
        A & 0 & 0 & 0 \\
        G & I & 0 & 0
    \end{bmatrix}
    \begin{bmatrix}
        \Delta x \\
        \Delta s \\
        \Delta y \\
        \Delta z
    \end{bmatrix}
    =
    \begin{bmatrix}
        -r_x \\
        -\mathrm{Arw}(s_k)^{-1}r_s \\
        -r_y \\
        -r_z
    \end{bmatrix}.
\end{equation}

However, the coefficient matrix above will only be symmetric if $\mathrm{Arw}(s_k)^{-1}\mathrm{Arw}(z_k)$ is symmetric. This can only be ensured if $\mc{K}$ only consisted of the non-negative orthant and no second-order cones. In other words, the system in \eqref{eq:asymmetric-system} would only be symmetric if Problem \ref{eq:problem} was a quadratic program rather than a second-order cone program. Even if we were to eliminate $\Delta s$ from Equation \eqref{eq:asymmetric-system}, we would still not have a symmetric linear system.
\section{Custom LDL factorization performance}\label{appendix:custom-ldl-perf}
This section provides empirical evidence supporting our claims in Sections \ref{sec:intro} and \ref{subsec:custom-linalg}. Specifically, we show that sparse linear algebra incurs extra overhead from identifying and locating nonzero elements, leading to more CPU instructions, memory accesses and cache misses. In contrast, a custom implementation that hardcodes the exact memory accesses and floating point operations reduces these inefficiencies. We specifically compare \textsc{qdldl} (used in {\qoco}) against our custom $LDL^\top$ factorization (used in {\qococ}), since the most expensive operation in a generic implementation of Algorithm \ref{alg:pdipm} is the $LDL^\top$ factorization. To compare the two matrix factorization implementations, we use Linux \texttt{perf}, a profiling tool that provides CPU performance metrics, which will allow us to quantify the computational overhead of sparse versus custom factorizations.

As a representative example of a matrix we would factor, we consider the KKT matrix associated with the Linear Quadratic Regulator (LQR). The LQR problem is

\begin{equation*}
    \begin{split}
        \underset{x, u}{\mathrm{minimize}} 
        \quad & \frac{1}{2}\left(\sum_{k = 1}^{T} x_k^\top Q_k x_k + \sum_{k = 1}^{T-1}u_k^\top R_k u_k\right) \\
        \mathrm{subject\;to} 
        \quad & x_{k+1} = A_kx_k + B_ku_k \;  \quad \forall k \in [1, T-1] \\
        & x_1 = x_{\mathrm{init}}, \\
    \end{split}
\end{equation*}

and its associated KKT matrix is 

\begin{equation*}
    K = \begin{bmatrix}
        P & H^\top \\
        H & -\epsilon I
    \end{bmatrix},
\end{equation*}

where

\begin{equation*}
    \begin{split}
        P &= \mathrm{blkdiag}(Q_1, Q_2, \ldots, Q_T, R_1, R_2, \ldots, R_{T-1}) \\
        H &= \begin{bmatrix}
            A_1 & -I & & & & &  B_1 \\
            & A_2 & -I & & & & &  B_2  \\
            & & & \ddots & & & & & \ddots \\
            & & & & A_{N-1} & -I & & & & B_{N-1} \\
            I & & & & & & & & & \\
            \end{bmatrix}.
    \end{split}
\end{equation*}
    
\begin{figure}[H]
    \captionsetup{labelfont=bf}
    \centering
    \includegraphics[width=0.75\textwidth]{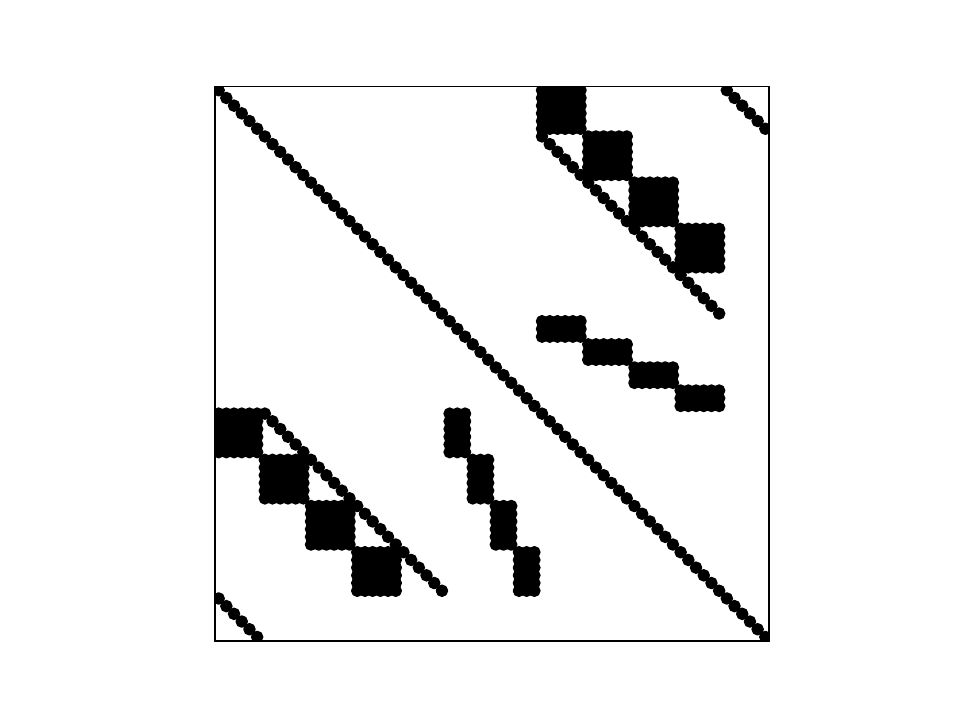}
    \caption{\bf Sparsity pattern of LQR KKT matrix for T=5}
    \label{fig:lqr-spy}
\end{figure}

Notice that we apply static regularization to ensure that the $(2,2)$ block of K is negative definite, making $K$ quasidefinite. The sparsity pattern of $K$ is depicted by Figure \ref{fig:lqr-spy}.

For simplicity, we assume time-invariant system and cost matrices, i.e., $Q_k=R_k=I$ and fixed matrices $A \in \R^{6 \times 6}$, $B \in \R^{6 \times 3}$. We choose the elements of $A$ and $B$ randomly and consider time horizons $T$: $\{5,15,50,75,100\}$. The actual elements of these matrices do not matter much, as our goal with these experiments is to understand the CPU behavior which (to first order) is unaffected by the value of elements in the matrices. To get accurate profiling results, we run each matrix factorization 1 million times to get the average runtime, number of instructions, number of L1 data cache loads, number of L1 cache misses, number of branches, and number of branch mispredictions.

From Table \ref{table:perf}, we can see that the custom $LDL^\top$ factorization is about four times faster than the sparse factorization because it does much less work: issuing fewer instructions, performing fewer memory accesses, and incurring fewer L1 cache misses. Although $\textsc{qdldl}$ has a far more branches than our custom $LDL^\top$ implementation, the branches are quite predictable, and pipeline stalls are avoided.

By hardcoding the exact memory accesses and floating point operations, our custom $LDL^\top$ factorization eliminates much of the overhead associated with sparse matrix operations. This substantial reduction in memory accesses and instructions leads to a significantly faster factorization.

\newpage
{
\scriptsize
\begin{longtable}{l||cc||cc||cc||cc||cc||cc||}
\captionsetup{labelfont=bf}
\caption{\bf \texttt{Perf} results for qdldl and custom $LDL^\top$ factorization \label{table:perf}} \\ 
\scriptsize
 & \multicolumn{2}{c||}{\underline{Runtime (us)}} & \multicolumn{2}{c||}{\underline{Instructions}} &\multicolumn{2}{c||}{\underline{L1 cache loads}} & \multicolumn{2}{c||}{\underline{L1 cache misses}} & \multicolumn{2}{c||}{\underline{Branches}} & \multicolumn{2}{c||}{\underline{Branch misses}}  \\ 
Size & qdldl & custom & qdldl & custom & qdldl & custom & qdldl & custom & qdldl & custom & qdldl & custom\\[1ex]
\hline
\endhead
5 & 2.103 & 0.435 \winner & 36650 & 4858 \winner & 14049 & 1851 \winner & 0 \winner & 0 \winner & 5419 & 5 \winner & 1 & 0 \winner \\ 
15 & 7.507 & 1.777 \winner & 132796 & 18053 \winner & 50550 & 6442 \winner & 572 & 27 \winner & 19124 & 8 \winner & 4 & 0 \winner \\ 
50 & 27.077 & 6.117 \winner & 469289 & 62046 \winner & 179709 & 21468 \winner & 2539 & 1204 \winner & 67081 & 17 \winner & 11 & 1 \winner \\ 
75 & 38.848 & 8.485 \winner & 709731 & 93912 \winner & 279033 & 32211 \winner & 3735 & 1828 \winner & 101347 & 24 \winner & 14 & 2 \winner \\ 
100 & 49.577 & 12.894 \winner & 950101 & 125747 \winner & 365981 & 43559 \winner & 5567 & 2459 \winner & 135597 & 35 \winner & 18 & 5 \winner \\ 
\end{longtable}

}
\section{Problem classes} \label{appendix:problem-classes}

In this section, we describe the five problem classes considered in our numerical experiments, along with various problem sizes and instances within each problem class. We use the notation $\mc{U}(a,b)$ to denote a uniform distribution on the interval $[a,b]$, $\mc{N}(\mu, \sigma)$ for a Gaussian distribution with mean $\mu$ and standard deviation $\sigma$, $0_{m \times n}$ for the matrix of all zeros in $\R^{m \times n}$, and $I_n$ for the identity matrix in $\R^{n \times n}$.

\subsection{Robust Kalman filter}
We consider the robust Kalman filtering problem for vehicle tracking, as outlined in \cite{cvxpyRobustKalmanFilter}. This problem is a quadratic objective second-order cone program (SOCP) and takes the form

\begin{equation}\label{eq:rkf}
\begin{split}
    \underset{x_k, w_k, v_k}{\text{minimize}} 
    \quad & \sum_{k=0}^{N-1}\left(\|w_k\|_2^2 + \tau \phi_{\rho}(v_k)\right) \\
    \text{subject to} 
    \quad & x_{k+1} = Ax_k + Bw_k \quad \forall k \in [0, N-1] \\
    \quad & y_{k} = Cx_k + v_k \quad \forall k \in [0, N-1], \\
\end{split}
\end{equation}

where $\phi_{\rho}$ is the Huber loss function

\begin{equation*}
    \phi_{\rho}(z) = 
    \begin{cases} 
      \|z\|_2^2 & \|z\|_2 \leq \rho \\
      2\rho\|z\|_2 - \rho^2 & \|z\|_2 > \rho. \\
   \end{cases}
\end{equation*}

In this problem, $x_k \in \R^4$ represents the state of the vehicle at timestep $k$,  $w_k \in \R^2$ is an unknown force vector acting on the vehicle at timestep $k$, $v_k \in \R^2$ is the measurement noise vector at timestep $k$, and $y_k \in \R^2$ represents a noisy measurement of the vehicle's state at timestep $k$. The matrices $A$ and $B$ are known dynamics matrices for the system, and $C$ is the known observation matrix. To reconstruct the vehicle's state trajectory, we solve Problem \ref{eq:rkf}. 

We use the system matrices

\begin{align*}
    A &= \begin{bmatrix}
        1 & 0 & \left(1-\frac{\gamma}{2}\Delta t\right)\Delta t & 0 \\
        0 & 1 & 0 & \left(1-\frac{\gamma}{2}\Delta t\right)\Delta t \\
        0 & 0 & 1 - \gamma \Delta t & 0 \\
        0 & 0 & 0 & 1 - \gamma \Delta t
    \end{bmatrix} \\
    \vspace{10mm}
    B &= \begin{bmatrix}
        \frac{1}{2} \Delta t^2 & 0 & \\
        0 & \frac{1}{2} \Delta t^2 \\
        \Delta t & 0 \\
        0 &  \Delta t
    \end{bmatrix} \\
    C &= \begin{bmatrix}
        1 & 0 & 0 & 0 \\
        0 & 1 & 0 & 0 \\
    \end{bmatrix}, \\
\end{align*}

which correspond to a two-dimensional double integrator model with linear velocity drag and unit mass, where the input is a force and we observe the position of the vehicle. We use $\gamma=0.05$ for the velocity drag parameter, $\Delta t$ for the discretization time, and $\rho=2$ and $\tau=2$ for the Huber loss function.

\paragraph{Problem sizes} To generate different problem sizes, we vary the number of timesteps $N$. We consider the following values for $N$: $\{25, 50, 75, 125, 175, 225, 300, 375, 450, 500\}$. For each of these problem sizes, we set $\Delta t = T/(N - 1)$, where $T=50$ is the time horizon.

\paragraph{Problem instances} For each problem size, we generate 20 random problem instances, corresponding to 20 different observation trajectories $y_t$. To do this, we first randomly generate the force and noise vectors $w_t$ and $v_t$, respectively, for  $t \in [0, N-1]$. Then, we forward-simulate an initial state of $x_0 = [0 \; 0 \; 0 \; 0]$ through the dynamics defined in the constraints of Problem \ref{eq:rkf}. Each element of $w_t$ and $v_t$ is drawn from a standard Gaussian distribution $\mc{N}(0, 1)$. Additionally, for 20 \% of the elements of $v_t$, we introduce outlier noise by drawing elements from $\mc{N}(0, 20)$, to simulate noise with a larger magnitude. 

\subsection{Lossless convexification}
The losslessly convexified powered-descent guidance problem is an optimal control problem that aims to compute a soft landing trajectory for a rocket while respecting relevant constraints. One important constraint is the lower thrust bound on the engine, which is nonconvex. However, it has been shown that this constraint admits a \textit{lossless} convex relaxation \cite{acikmese2007lcvx, acikmese2013lossless}. This means that the solution to the relaxed convex problem guarantees a solution to the original nonconvex problem. The losslessly convexified problem can be written as an SOCP

\begin{equation*}
    \begin{split}
        \underset{x, z, u, \sigma}{\text{minimize}} 
        \quad & -z_T  \\
        \text{subject to} 
        \quad & x_{k+1} = Ax_k + Bu_k + g \quad \forall k \in [0, T-1] \\  
        \quad & z_{k+1} = z_k - \alpha \sigma_k \Delta t \quad \forall k \in [0, T-1] \\  
        \quad & \|u_k\|_2 \leq \sigma_k \quad \forall k \in [0, T-1] \\
        \quad & \log(m_{\mathrm{wet}} - \alpha\rho_2 k \Delta t) \leq z_k \leq \log(m_{\mathrm{wet}} - \alpha\rho_1 k \Delta t) \quad \forall k \in [0, T-1] \\
        \quad & \mu_{1,k}\left[1-[z_k-z_{0,k}] + \frac{[z_k-z_{0,k}]^2}{2}\right] \leq \sigma_k \leq \mu_{2,k}[1-(z_k-z_{0,k})] \quad \forall k \in [0, T-1] \\ 
        \quad & e_3^\top u_k \geq \sigma_k \cos(\theta_{\mathrm{max}}) \quad \forall k \in [0, T-1] \\
        \quad & x_0 = x_{\mathrm{init}}, \; z_0 = \log(m_{\mathrm{wet}}), \; z_T \geq \log(m_{\mathrm{dry}}),
    \end{split}
\end{equation*}

where $x_k \in \R^6$ is the position and velocity of the rocket at timestep $k$, $u_k \in \R^3$ is the thrust per unit mass at timestep $k$, $z_k \in \R$ is the natural logarithm of the rocket's mass at timestep $k$, and $\sigma_k \in \R$ is a slack variable at timestep $k$. Further we define $z_{0,k} \coloneqq \log(m_{\mathrm{wet}} - \alpha \rho_2 k \Delta t)$, $\mu_{1,k} \coloneqq \rho_1 {\rm e}^{-z_0}$, and $\mu_{2,k} \coloneqq \rho_2 {\rm e}^{-z_0}$. The parameter $\rho_1 \in \R$ is the minimum thrust of the engine, $\rho_2 \in \R$ is the maximum thrust of the engine, $\alpha \in \R$ is a mass depletion parameter that describes how efficient the engine is, $m_{\mathrm{wet}} \in \R$ is the sum of the dry mass of the rocket and the initial fuel mass, $\theta_{\mathrm{max}} \in \R$ is the maximum angle from vertical that the thrust vector can make, $x_{\mathrm{init}} \in \R^6$ is the initial state of the rocket, $\Delta t \in \R$ is the discretization time, and $e_3$ is the third canonical basis vector in $\R^3$. The dynamics matrices $A$ and $B$, and the vector $g$, are defined as

\begin{align*}
    A &= \begin{bmatrix}
        I_3 & \Delta t (I_3) \\
        0_{3 \times 3} & I_3 
    \end{bmatrix} \\
    \vspace{10mm}
    B &= \begin{bmatrix}
        \frac{1}{2} \Delta t^2 (I_3) \\
        \Delta t (I_3)
    \end{bmatrix} \\
    g &= [0 \;\; 0 \;\; -0.5 g_0 \Delta t ^2 \;\; 0 \;\; 0 \;\; -g_0 \Delta t].
\end{align*}

For this problem class we use the parameters

\begin{equation*}
    g_0 = 9.807, \; \rho_1 = 100, \; \rho _2 = 500, \; m_{\mathrm{dry}} = 25, \; m_{\mathrm{wet}} = 35, \; \theta_{\mathrm{max}} = \pi / 4, \; \alpha = 0.001.  \\
\end{equation*}

\paragraph{Problem sizes} To generate different problem sizes, we vary the number of timesteps $T$. The discretization time is computed as $\Delta t = t_f / (T - 1)$, where $t_f = 20$. We consider the following values for $T$: $\{15, 50, 75, 100, 125, 150, 200, 250, 300, 350\}$.

\paragraph{Problem instances} For each problem size (i.e., for each number of timesteps), we generate 20 random problem instances, which involves generating a random initial state $x_{\mathrm{init}}$. We generate $x_{\mathrm{init}}$ as

\begin{equation*}
    x_{\mathrm{init}} \sim \begin{bmatrix} \mc{U}(-10, 10) \\ \mc{U}(-10, 10) \\ \mc{U}(200, 400) \\ 0 \\ 0 \\ 0 \end{bmatrix}.
\end{equation*}

\subsection{Group lasso}
The group lasso problem is a quadratic objective SOCP. It is a variation of the least-squares regression problem, where the regression variables are partitioned into disjoint groups, and a regularization term is added to encourage group sparsity. Let $x = [x^{(1)}, x^{(2)}, \ldots, x^{(N)}]$ represent the partitioning of the regression variables, where each $x^{(i)}$ is a group of variables. We can write the problem as

\begin{equation*}
    \begin{split}
        \underset{x}{\text{minimize}} 
        \quad & \|Ax - b\|_2^2 + \lambda \sum_{i=1}^{N} \|x^{(i)}\|_2.  \\
    \end{split}
\end{equation*}

This can be written as an SOCP by introducing a slack variable $t \in \R^N$, and an auxiliary variable $y \in \R^m$ to avoid computing $A^\top A$ and maintain sparsity in the Hessian of the objective function. The problem is now  

\begin{equation*}
    \begin{split}
        \underset{x, t}{\text{minimize}} 
        \quad & \|y\|_2^2 + \lambda \sum_{i=1}^{N} t_i  \\
        \text{subject to} 
        \quad & \|x^{(i)}\|_2 \leq t_i \quad \forall i \in [1, N] \\  
        \quad & y = Ax - b,
    \end{split}
\end{equation*}

where $A \in \R^{m \times n}$, $b \in \R^m$.

\paragraph{Problem sizes} To generate different problem sizes, we vary the number of groups $N$. We consider the following values for $N$: $\{1, 2, 3, 4, 5, 8, 10, 12, 14, 16\}$.

\paragraph{Problem instances} For each problem size (i.e., for each number of groups) we generate 20 random problem instances, corresponding to the generation of the matrix $A$ and the vector $b$. We set $A$ to have $m = 250N$ rows and $n = 10N$ columns, where $x$ is partitioned into $N$ groups, each containing ten regression variables. The matrix $A$ is sparse, with 10\% of its elements nonzero, and the nonzero elements are drawn from the uniform distribution $\mc{U}(0, 1)$.

To generate $b$, we first generate $\hat{x} \in \R^n$ and partition it into $N$ groups of ten variables. Half of the groups are set to zero, and the remaining elements are drawn from the standard Gaussian distribution $\mc{N}(0, 1)$. Next, we generate an error vector $e \in \R^m$ with elements drawn from $\mc{N}(0, 1/n)$. Finally, we compute $b = A\hat{x} + e$. Note that each problem instance will have a different sparsity pattern for $A$, which requires us to generate a new custom solver for each instance.

\subsection{Portfolio optimization}
We consider the Markowitz portfolio optimization problem with a no-short-selling constraint \cite{markowitz1952-ds,Boyd2004}. This problem is the QP

\begin{equation*}
    \begin{split}
        \underset{x}{\text{minimize}} 
        \quad & \gamma x^\top \Sigma x - \mu^\top x  \\
        \text{subject to} 
        \quad & \sum_{i=1}^{n} x_i = 1 \\  
        \quad &  x_i \geq 0 \quad \forall i \in [1, n],\\    
    \end{split}
\end{equation*}

where $x \in \R^n$ is the allocation vector, with $x_i$ representing the percentage of resources invested in asset $i$, $\mu \in \R^n$ is the vector of expected returns, and $\Sigma \in \mb{S}^n_+$ is the covariance matrix of asset returns. 

Typically, the return covariance matrix $\Sigma$ is approximated as the sum of a low-rank matrix and a diagonal matrix. Thus, we write $\Sigma = FF^\top + D$, where $F \in \R^{n \times k}$ is a factor matrix of rank $k$. We then add the auxiliary variable $y$ and the constraint $y = F^\top x$ to avoid computing $FF^\top$ and maintain sparsity in the Hessian of the objective function. The final problem can be rewritten as

\begin{equation}\label{eq:portfolio}
    \begin{split}
        \underset{x, y}{\text{minimize}} 
        \quad & x^\top D x + y^\top y - \gamma^{-1}\mu^\top x  \\
        \text{subject to} 
        \quad & y = F^\top x \\  
        \quad & \sum_{i=1}^{n} x_i = 1 \\  
        \quad &  x_i \geq 0 \quad \forall i \in [1, n].\\    
    \end{split}
\end{equation}

\paragraph{Problem sizes} To generate different problem sizes, we vary the number of factors, $k$. We consider the following values for $k$: $\{2, 4, 6, 8, 10, 15, 20, 25, 30, 35\}$.

\paragraph{Problem instances} For each problem size (i.e., each number of factors), we generate 20 random problem instances. This corresponds to generating the matrices $D$, $F$, and the vector $\mu$. We set the number of assets to $n = 100k$, and for all instances, we take $\gamma = 1$. The diagonal elements of $D$ are drawn from $\mc{N}(0, 1)$, and the elements of $\mu$ are drawn from $\mc{N}(0, \sqrt{k})$. The matrix $F$ is a sparse matrix, with 50\% of its elements being nonzero, and its nonzero elements are drawn from $\mc{N}(0, 1)$. Note that each problem instance will have a different sparsity pattern for $F$, which requires us to generate a new custom solver for each instance.

\subsection{Oscillating masses}
The oscillating masses problem is an optimal control problem involving a system of $N$ masses in one dimension, where each mass is connected to its neighboring mass by a spring, and the outermost two masses are attached to a fixed wall. The goal is to apply forces to each mass to bring the system to equilibrium \cite{wang2010fast}. This problem can be formulated as a QP

\begin{equation*}
    \begin{split}
        \underset{x, u}{\mathrm{minimize}} 
        \quad & \frac{1}{2}\left(\sum_{k = 0}^{T} x_k^\top Q x_k + \sum_{k = 0}^{T-1}u_k^\top R u_k\right) \\
        \mathrm{subject\;to} 
        \quad & x_{k+1} = Ax_k + Bu_k \;  \quad \forall k \in [0, T-1] \\
        & x_0 = x_{\mathrm{init}} \\
        & \|x_k\|_\infty \leq x_{\mathrm{max}} \quad \forall k \in [0, T] \\
        & \|u_k\|_\infty \leq u_{\mathrm{max}} \quad \forall k \in [0, T-1].
    \end{split}
\end{equation*}

where $x_k \in \R^{2N}$ is the position and velocity of the masses at timestep $k$, $u_k \in \R^N$ is force applied to each mass at timestep $k$, and $Q \in \mb{S}^{2N}_+$ and $R \in \mb{S}^{N}_+$ penalize deviation from equilibrium and control effort respectively.

The dynamics matrices $A$ and $B$ are derived from an exact discretization of the continuous-time linear dynamics, using a zero-order hold on control input

\begin{align*}
    A &= {\rm e}^{A_c \Delta t} \\
    \vspace{10mm}
    B &= A_c^{-1}(A - I_{2N})B_c, \\
\end{align*}

where

\begin{align*}
    A_c &= \begin{bmatrix}
        0_{N \times N} & I_N \\
        L_N & 0_{N \times N} 
    \end{bmatrix} \\
    \vspace{10mm}
    B_c &= \begin{bmatrix}
        0_{N \times N} \\
        I_N
    \end{bmatrix}. \\
\end{align*}

Here, $L_N \in \R^{N \times N}$ is a symmetric tridiagonal matrix with $-2$ on the main diagonal and $1$ on the subdiagonal and superdiagonal.

For this problem class, we use the parameters

\begin{equation*}
    N = 4, \; \Delta t = 0.25, \; x_{\mathrm{max}} = 2, \; u_{\mathrm{max}} = 5.
\end{equation*}

\paragraph{Problem sizes} To generate different problem sizes, we vary the number of timesteps $T$. We consider the following values for $T$: $\{8, 20, 32, 44, 56, 76, 96, 116, 136, 156\}$.

\paragraph{Problem instances} For each problem size (i.e., each number of timesteps), we generate 20 random problem instances. This involves generating a random initial state $x_{\mathrm{init}}$, as well as the matrices $Q$ and $R$. We draw each element of $x_{\mathrm{init}}$ from $\mc{N}(0, 1)$, then clip the elements such that they lie within the interval $[-0.9x_{\mathrm{max}}, 0.9x_{\mathrm{max}}]$. This ensures that the initial state satisfies the state constraints, making the resulting problem feasible. The matrices $Q$ and $R$ are generated as diagonal matrices, with each diagonal element drawn from the uniform distribution $\mc{U}(0, 10)$.

\section{Detailed numeric results}\label{appendix:detailed-results}
\newpage

\begin{landscape}
    \scriptsize
\begin{longtable}{lc||cccccc||cccccc||}
\captionsetup{labelfont=bf}
\caption{\bf Iterations and solver runtimes for robust Kalman filter problems} \\ 
 & &  \multicolumn{6}{c||}{\underline{Iterations}} & \multicolumn{6}{c||}{\underline{Solver Runtime (s)}}\\[2ex] 
Problem & Size & \textsc{clarabel} & \textsc{ecos} & \textsc{gurobi} & \textsc{mosek} & \textsc{qoco} & \textsc{qoco\_custom} & \textsc{clarabel} & \textsc{ecos} & \textsc{gurobi} & \textsc{mosek} & \textsc{qoco} & \textsc{qoco\_custom} \\[1ex]
\hline
\endhead
\textsc{robust\_kalman\_filter\_N\_25\_i\_0} & 775 & 12 & 15 & -&  \winner 9 &  \winner 9 &  \winner 9 & 0.00082 & 0.00148 & 0.00618 & 0.00123 & 0.00043 &  \winner 0.00018 \\ 
\textsc{robust\_kalman\_filter\_N\_25\_i\_1} & 775 & 10 & 15 & -& 10 &  \winner 9 &  \winner 9 & 0.00071 & 0.00154 & 0.00652 & 0.00129 & 0.00046 &  \winner 0.00018 \\ 
\textsc{robust\_kalman\_filter\_N\_25\_i\_2} & 775 & 11 & 18 & -& 11 &  \winner 10 &  \winner 10 & 0.00077 & 0.00175 & 0.00701 & 0.00134 & 0.00047 &  \winner 0.00020 \\ 
\textsc{robust\_kalman\_filter\_N\_25\_i\_3} & 775 & 10 & 17 & -& 10 &  \winner 9 &  \winner 9 & 0.00072 & 0.00171 & 0.00616 & 0.00129 & 0.00045 &  \winner 0.00018 \\ 
\textsc{robust\_kalman\_filter\_N\_25\_i\_4} & 775 & 11 & 15 & -& 11 &  \winner 9 &  \winner 9 & 0.00077 & 0.00148 & 0.00630 & 0.00137 & 0.00043 &  \winner 0.00018 \\ 
\textsc{robust\_kalman\_filter\_N\_25\_i\_5} & 775 & 10 & 14 & -& 9 &  \winner 8 &  \winner 8 & 0.00071 & 0.00145 & 0.00659 & 0.00124 & 0.00041 &  \winner 0.00016 \\ 
\textsc{robust\_kalman\_filter\_N\_25\_i\_6} & 775 & 11 & 15 & -& 10 &  \winner 9 &  \winner 9 & 0.00076 & 0.00152 & 0.00656 & 0.00130 & 0.00045 &  \winner 0.00018 \\ 
\textsc{robust\_kalman\_filter\_N\_25\_i\_7} & 775 & 10 & 17 & -& 11 &  \winner 9 &  \winner 9 & 0.00071 & 0.00168 & 0.00728 & 0.00136 & 0.00045 &  \winner 0.00018 \\
\end{longtable}

    \scriptsize
\begin{longtable}{lc||cccccc||cccccc||}
\captionsetup{labelfont=bf}
\caption{\bf Iterations and solver runtimes for lossless convexification problems} \\ 
 & &  \multicolumn{6}{c||}{\underline{Iterations}} & \multicolumn{6}{c||}{\underline{Solver Runtime (s)}}\\[2ex] 
Problem & Size & \textsc{clarabel} & \textsc{ecos} & \textsc{gurobi} & \textsc{mosek} & \textsc{qoco} & \textsc{qoco\_custom} & \textsc{clarabel} & \textsc{ecos} & \textsc{gurobi} & \textsc{mosek} & \textsc{qoco} & \textsc{qoco\_custom} \\[1ex]
\hline
\endhead
\textsc{lcvx\_N\_15\_i\_0} & 644 & 8 &  \winner 7 & -& 9 & 9 & 9 & 0.00074 & 0.00060 & 0.00438 & 0.00111 & 0.00049 &  \winner 0.00023 \\ 
\textsc{lcvx\_N\_15\_i\_1} & 644 & 8 &  \winner 7 & -& 9 & 9 & 9 & 0.00072 & 0.00060 & 0.00428 & 0.00110 & 0.00049 &  \winner 0.00024 \\ 
\textsc{lcvx\_N\_15\_i\_2} & 644 & 10 &  \winner 7 & -& 9 & 10 & 10 & 0.00085 & 0.00062 & 0.00458 & 0.00099 & 0.00055 &  \winner 0.00025 \\ 
\textsc{lcvx\_N\_15\_i\_3} & 644 & 8 &  \winner 7 & -& 9 & 10 & 10 & 0.00073 & 0.00057 & 0.00442 & 0.00107 & 0.00052 &  \winner 0.00026 \\ 
\textsc{lcvx\_N\_15\_i\_4} & 644 & 9 &  \winner 7 & -& 9 & 9 & 9 & 0.00078 & 0.00062 & 0.00472 & 0.00109 & 0.00049 &  \winner 0.00023 \\ 
\textsc{lcvx\_N\_15\_i\_5} & 644 & 10 &  \winner 7 & -& 10 & 9 & 9 & 0.00084 & 0.00060 & 0.00446 & 0.00115 & 0.00049 &  \winner 0.00023 \\ 
\textsc{lcvx\_N\_15\_i\_6} & 644 &  \winner 7 &  \winner 7 & -& 9 & 10 & 10 & 0.00068 & 0.00058 & 0.00417 & 0.00109 & 0.00053 &  \winner 0.00025 \\ 
\textsc{lcvx\_N\_15\_i\_7} & 644 & 9 &  \winner 7 & -& 9 & 10 & 10 & 0.00078 & 0.00060 & 0.00450 & 0.00108 & 0.00053 &  \winner 0.00025 \\
\end{longtable}

    \scriptsize
\begin{longtable}{lc||cccccc||cccccc||}
\captionsetup{labelfont=bf}
\caption{\bf Iterations and solver runtimes for group lasso problems} \\ 
 & &  \multicolumn{6}{c||}{\underline{Iterations}} & \multicolumn{6}{c||}{\underline{Solver Runtime (s)}}\\[2ex] 
Problem & Size & \textsc{clarabel} & \textsc{ecos} & \textsc{gurobi} & \textsc{mosek} & \textsc{qoco} & \textsc{qoco\_custom} & \textsc{clarabel} & \textsc{ecos} & \textsc{gurobi} & \textsc{mosek} & \textsc{qoco} & \textsc{qoco\_custom} \\[1ex]
\hline
\endhead
\textsc{group\_lasso\_N\_1\_i\_0} & 761 &  \winner 4 & 11 & -&  \winner 4 & 6 & 6 & 0.00025 & 0.00111 & 0.00564 & 0.00047 & 0.00020 &  \winner 0.00007 \\ 
\textsc{group\_lasso\_N\_1\_i\_1} & 761 &  \winner 4 & 11 & -& 5 & 6 & 6 & 0.00027 & 0.00109 & 0.00580 & 0.00049 & 0.00020 &  \winner 0.00007 \\ 
\textsc{group\_lasso\_N\_1\_i\_2} & 761 &  \winner 4 & 10 & -& 5 & 6 & 6 & 0.00027 & 0.00098 & 0.00597 & 0.00049 & 0.00021 &  \winner 0.00007 \\ 
\textsc{group\_lasso\_N\_1\_i\_3} & 761 &  \winner 4 & 12 & -& 5 & 6 & 6 & 0.00027 & 0.00122 & 0.00577 & 0.00050 & 0.00020 &  \winner 0.00007 \\ 
\textsc{group\_lasso\_N\_1\_i\_4} & 761 &  \winner 4 & 11 & -& 5 & 6 & 6 & 0.00027 & 0.00110 & 0.00545 & 0.00050 & 0.00020 &  \winner 0.00007 \\ 
\textsc{group\_lasso\_N\_1\_i\_5} & 761 &  \winner 4 & 11 & -& 5 & 6 & 6 & 0.00026 & 0.00113 & 0.00610 & 0.00048 & 0.00020 &  \winner 0.00007 \\ 
\textsc{group\_lasso\_N\_1\_i\_6} & 761 &  \winner 4 & 11 & -& 5 & 6 & 6 & 0.00027 & 0.00109 & 0.00588 & 0.00049 & 0.00022 &  \winner 0.00006 \\ 
\textsc{group\_lasso\_N\_1\_i\_7} & 761 &  \winner 5 & 10 & -&  \winner 5 & 6 & 6 & 0.00030 & 0.00099 & 0.00549 & 0.00049 & 0.00021 &  \winner 0.00007 \\
\end{longtable}

    \scriptsize
\begin{longtable}{lc||cccccc||cccccc||}
\captionsetup{labelfont=bf}
\caption{\bf Iterations and solver runtimes for portfolio optimization problems} \\ 
 & &  \multicolumn{6}{c||}{\underline{Iterations}} & \multicolumn{6}{c||}{\underline{Solver Runtime (s)}}\\[2ex] 
Problem & Size & \textsc{clarabel} & \textsc{ecos} & \textsc{gurobi} & \textsc{mosek} & \textsc{qoco} & \textsc{qoco\_custom} & \textsc{clarabel} & \textsc{ecos} & \textsc{gurobi} & \textsc{mosek} & \textsc{qoco} & \textsc{qoco\_custom} \\[1ex]
\hline
\endhead
\textsc{portfolio\_N\_2\_i\_0} & 804 &  \winner 9 & 12 & 14 & 11 &  \winner 9 &  \winner 9 & 0.00041 & 0.00087 & 0.00205 & 0.00100 & 0.00025 &  \winner 0.00012 \\ 
\textsc{portfolio\_N\_2\_i\_1} & 804 & 9 & 14 & 13 & 11 &  \winner 8 &  \winner 8 & 0.00041 & 0.00096 & 0.00200 & 0.00101 & 0.00024 &  \winner 0.00010 \\ 
\textsc{portfolio\_N\_2\_i\_2} & 804 &  \winner 8 & 12 & 12 & 11 &  \winner 8 &  \winner 8 & 0.00038 & 0.00089 & 0.00204 & 0.00100 & 0.00024 &  \winner 0.00010 \\ 
\textsc{portfolio\_N\_2\_i\_3} & 804 &  \winner 8 & 15 & 16 & 11 &  \winner 8 &  \winner 8 & 0.00038 & 0.00110 & 0.00203 & 0.00101 & 0.00024 &  \winner 0.00011 \\ 
\textsc{portfolio\_N\_2\_i\_4} & 804 & 9 & 14 & 15 & 9 &  \winner 8 &  \winner 8 & 0.00041 & 0.00099 & 0.00207 & 0.00092 & 0.00024 &  \winner 0.00011 \\ 
\textsc{portfolio\_N\_2\_i\_5} & 804 & 10 & 14 & 14 & 13 &  \winner 9 &  \winner 9 & 0.00046 & 0.00097 & 0.00203 & 0.00110 & 0.00026 &  \winner 0.00011 \\ 
\textsc{portfolio\_N\_2\_i\_6} & 804 &  \winner 10 & 18 & 11 &  \winner 10 &  \winner 10 &  \winner 10 & 0.00045 & 0.00133 & 0.00198 & 0.00096 & 0.00028 &  \winner 0.00013 \\ 
\textsc{portfolio\_N\_2\_i\_7} & 804 & 8 & 10 & 11 & 9 &  \winner 7 &  \winner 7 & 0.00037 & 0.00078 & 0.00202 & 0.00091 & 0.00021 &  \winner 0.00009 \\
\end{longtable}

    \scriptsize
\begin{longtable}{lc||ccccccc||ccccccc||}
\captionsetup{labelfont=bf}
\caption{\bf Iterations and solver runtimes for oscillating masses problems} \\ 
 & &  \multicolumn{7}{c||}{\underline{Iterations}} & \multicolumn{7}{c||}{\underline{Solver Runtime (s)}}\\[2ex] 
Problem & Size & \textsc{clarabel} & \textsc{ecos} & \textsc{gurobi} & \textsc{mosek} & \textsc{qoco} & \textsc{qoco\_custom} & \textsc{cvxgen} & \textsc{clarabel} & \textsc{ecos} & \textsc{gurobi} & \textsc{mosek} & \textsc{qoco} & \textsc{qoco\_custom} & \textsc{cvxgen} \\[1ex]
\hline
\endhead
\textsc{oscillating\_masses\_N\_8\_i\_0} & 1344 &  \winner 5 & 13 & 11 & 11 &  \winner 5 &  \winner 5 & 7 & 0.00049 & 0.00108 & 0.00222 & 0.00187 & 0.00027 &  \winner 0.00011 & 0.00014 \\ 
\textsc{oscillating\_masses\_N\_8\_i\_1} & 1344 &  \winner 5 & 14 & 8 & 10 &  \winner 5 &  \winner 5 & 6 & 0.00048 & 0.00117 & 0.00218 & 0.00181 & 0.00028 &  \winner 0.00012 & 0.00013 \\ 
\textsc{oscillating\_masses\_N\_8\_i\_2} & 1344 &  \winner 5 & 17 & 9 & 11 &  \winner 5 &  \winner 5 & 6 & 0.00048 & 0.00139 & 0.00216 & 0.00189 & 0.00026 &  \winner 0.00011 & 0.00013 \\ 
\textsc{oscillating\_masses\_N\_8\_i\_3} & 1344 &  \winner 5 & 14 & 8 & 11 &  \winner 5 &  \winner 5 & 6 & 0.00047 & 0.00115 & 0.00217 & 0.00190 & 0.00026 &  \winner 0.00011 & 0.00012 \\ 
\textsc{oscillating\_masses\_N\_8\_i\_4} & 1344 &  \winner 5 & 12 & 10 & 11 &  \winner 5 &  \winner 5 & 7 & 0.00049 & 0.00099 & 0.00219 & 0.00189 & 0.00029 &  \winner 0.00011 & 0.00015 \\ 
\textsc{oscillating\_masses\_N\_8\_i\_5} & 1344 &  \winner 5 & 12 & 8 & 11 &  \winner 5 &  \winner 5 & 6 & 0.00047 & 0.00100 & 0.00216 & 0.00190 & 0.00026 &  \winner 0.00011 & 0.00013 \\ 
\textsc{oscillating\_masses\_N\_8\_i\_6} & 1344 & 10 & 13 & 12 & 11 &  \winner 7 &  \winner 7 & 10 & 0.00083 & 0.00110 & 0.00223 & 0.00188 & 0.00036 &  \winner 0.00015 & 0.00020 \\ 
\textsc{oscillating\_masses\_N\_8\_i\_7} & 1344 &  \winner 5 & 13 & 9 & 10 &  \winner 5 &  \winner 5 & 6 & 0.00048 & 0.00104 & 0.00219 & 0.00178 & 0.00027 &  \winner 0.00011 & 0.00012 \\
\end{longtable}

    \scriptsize
\begin{longtable}{lc||cccccc||cccccc||}
\captionsetup{labelfont=bf}
\caption{\bf Iterations and solver runtimes for mpc problems} \\ 
 & & \multicolumn{6}{c||}{\underline{Iterations}} & \multicolumn{6}{c||}{\underline{Solver Runtime (s)}}\\[2ex]
Problem & Size & \textsc{qoco} & \textsc{qoco\_custom} & \textsc{clarabel} & \textsc{ecos} & \textsc{gurobi} & \textsc{mosek} & \textsc{qoco} & \textsc{qoco\_custom} & \textsc{clarabel} & \textsc{ecos} & \textsc{gurobi} & \textsc{mosek} \\[1ex]
\hline
\endhead
\textsc{aircraft\_1} & 504 &  \winner 8 &  \winner 8 & 9 & 15 & 10 & 11 & 0.00018 &  \winner 0.00007 & 0.00036 & 0.00053 & 0.00202 & 0.00074 \\
\textsc{aircraft\_10} & 584 &  \winner 8 &  \winner 8 & 9 & 15 & 10 & 13 & 0.00020 &  \winner 0.00008 & 0.00041 & 0.00055 & 0.00204 & 0.00081 \\
\textsc{aircraft\_11} & 636 & 11 & 11 & 12 & 23 &  \winner 9 & 14 & 0.00026 &  \winner 0.00010 & 0.00054 & 0.00097 & 0.00203 & 0.00099 \\
\textsc{aircraft\_12} & 636 & 11 & 11 & 12 & 23 &  \winner 9 & 14 & 0.00026 &  \winner 0.00011 & 0.00055 & 0.00093 & 0.00202 & 0.00100 \\
\textsc{aircraft\_13} & 3116 & 14 & 14 & 15 & 25 &  \winner 11 & 19 & 0.00141 &  \winner 0.00073 & 0.00334 & 0.00523 & 0.00276 & 0.00965 \\
\textsc{aircraft\_2} & 524 &  \winner 8 &  \winner 8 & 9 & 17 & 10 & 13 & 0.00017 &  \winner 0.00007 & 0.00035 & 0.00053 & 0.00202 & 0.00085 \\
\textsc{aircraft\_3} & 464 &  \winner 7 &  \winner 7 &  \winner 7 & 14 & 16 & 12 & 0.00015 &  \winner 0.00005 & 0.00025 & 0.00048 & 0.00209 & 0.00070 \\
\textsc{aircraft\_4} & 504 &  \winner 8 &  \winner 8 & 9 & 15 & 10 & 11 & 0.00018 &  \winner 0.00007 & 0.00035 & 0.00052 & 0.00203 & 0.00076 \\
\end{longtable}

    \scriptsize
\begin{longtable}{lc||cccccc||cccccc||}
\captionsetup{labelfont=bf}
\caption{\bf Iterations and solver runtimes for mpc problems on Raspberry Pi CM4} \\ 
 & & \multicolumn{6}{c||}{\underline{Iterations}} & \multicolumn{6}{c||}{\underline{Solver Runtime (s)}}\\[2ex]
Problem & Size & \textsc{qoco} & \textsc{qoco\_custom} & \textsc{clarabel} & \textsc{ecos} & \textsc{gurobi} & \textsc{mosek} & \textsc{qoco} & \textsc{qoco\_custom} & \textsc{clarabel} & \textsc{ecos} & \textsc{gurobi} & \textsc{mosek} \\[1ex]
\hline
\endhead
\textsc{aircraft\_1} & 504 &  \winner 8 &  \winner 8 & 9 & 15 & 10 & 11 & 0.00119 &  \winner 0.00062 & 0.00269 & 0.00412 & 0.00340 & 0.00707 \\
\textsc{aircraft\_10} & 584 &  \winner 8 &  \winner 8 & 9 & 15 & 10 & 13 & 0.00138 &  \winner 0.00077 & 0.00313 & 0.00446 & 0.00347 & 0.00794 \\
\textsc{aircraft\_11} & 636 & 11 & 11 & 12 & 23 &  \winner 9 & 14 & 0.00179 &  \winner 0.00109 & 0.00412 & 0.00736 & 0.00353 & 0.00910 \\
\textsc{aircraft\_12} & 636 & 11 & 11 & 12 & 23 &  \winner 9 & 14 & 0.00180 &  \winner 0.00106 & 0.00412 & 0.00736 & 0.00357 & 0.00911 \\
\textsc{aircraft\_13} & 3116 & 14 & 14 & 15 & 25 &  \winner 11 & 18 & 0.01000 &  \winner 0.00772 & 0.02586 & 0.04058 & 0.01110 & 0.07582 \\
\textsc{aircraft\_2} & 524 &  \winner 8 &  \winner 8 & 9 & 17 & 10 & 13 & 0.00116 &  \winner 0.00063 & 0.00270 & 0.00427 & 0.00332 & 0.00795 \\
\textsc{aircraft\_3} & 464 &  \winner 7 &  \winner 7 &  \winner 7 & 14 & 16 & 12 & 0.00094 &  \winner 0.00046 & 0.00195 & 0.00353 & 0.00376 & 0.00642 \\
\textsc{aircraft\_4} & 504 &  \winner 8 &  \winner 8 & 9 & 15 & 10 & 11 & 0.00116 &  \winner 0.00062 & 0.00266 & 0.00420 & 0.00342 & 0.00707 \\
\end{longtable}

    \scriptsize
\begin{longtable}{lc||ccccc||ccccc||}
\captionsetup{labelfont=bf}
\caption{\bf Iterations and solver runtimes for Maros–Mészáros problems} \\ 
 & & \multicolumn{5}{c||}{\underline{Iterations}} & \multicolumn{5}{c||}{\underline{Solver Runtime (s)}}\\[2ex]
Problem & Size & \textsc{qoco} & \textsc{clarabel} & \textsc{ecos} & \textsc{gurobi} & \textsc{mosek} & \textsc{qoco} & \textsc{clarabel} & \textsc{ecos} & \textsc{gurobi} & \textsc{mosek} \\[1ex]
\hline
\endhead
\textsc{AUG2D} & 59800 &  \winner 0 &  \winner 0 & -& 2 & 14 &  \winner 0.00958 & 0.01984 & -& 0.02289 & 0.31375 \\
\textsc{AUG2DC} & 60200 &  \winner 0 &  \winner 0 & -& 2 & 12 &  \winner 0.01347 & 0.01958 & -& 0.02400 & 0.23210 \\
\textsc{AUG2DCQP} & 80400 &  \winner 11 & 12 & -& 17 & -& 0.08231 & 0.13210 & -&  \winner 0.05561 & -\\
\textsc{AUG2DQP} & 80000 &  \winner 11 & 13 & -& 19 & -& 0.07272 & 0.13726 & -&  \winner 0.05924 & -\\
\textsc{AUG3D} & 9219 &  \winner 0 &  \winner 0 & 13 & 2 & 8 &  \winner 0.00148 & 0.00230 & 0.03595 & 0.00389 & 0.03126 \\
\textsc{AUG3DC} & 10419 &  \winner 0 &  \winner 0 & 14 & 2 & 7 &  \winner 0.00201 & 0.00302 & 0.03647 & 0.00928 & 0.01731 \\
\textsc{AUG3DCQP} & 14292 &  \winner 9 & 10 & 29 & 16 & 10 &  \winner 0.01278 & 0.01802 & 0.08718 & 0.01828 & 0.02658 \\
\textsc{AUG3DQP} & 13092 &  \winner 11 &  \winner 11 & 31 & 18 & 17 &  \winner 0.01532 & 0.01932 & 0.08252 & 0.01974 & 0.10652 \\
\end{longtable}

    \scriptsize
\begin{longtable}{lc||ccccc||ccccc||}
\captionsetup{labelfont=bf}
\caption{\bf Iterations and solver runtimes for SuiteSparse problems} \\ 
 & & \multicolumn{5}{c||}{\underline{Iterations}} & \multicolumn{5}{c||}{\underline{Solver Runtime (s)}}\\[2ex]
Problem & Size & \textsc{qoco} & \textsc{clarabel} & \textsc{ecos} & \textsc{gurobi} & \textsc{mosek} & \textsc{qoco} & \textsc{clarabel} & \textsc{ecos} & \textsc{gurobi} & \textsc{mosek} \\[1ex]
\hline
\endhead
\textsc{ANSYS\_Delor295K\_huber} & 4471461 & 7 & 7 & 8 & -&  \winner 4 &  \winner 3.12931 & 5.19758 & 10.03436 & -& 3.17687 \\
\textsc{ANSYS\_Delor295K\_lasso} & 10288503 &  \winner 7 &  \winner 7 & -& 12 & 33 &  \winner 7.50801 & 12.16215 & -& 12.36872 & 24.15621 \\
\textsc{ANSYS\_Delor338K\_huber} & 6614251 & 7 & 7 & -&  \winner 0 & 4 & 4.05127 & 6.12029 & -&  \winner 0.92594 & 3.42096 \\
\textsc{ANSYS\_Delor338K\_lasso} & 8446303 & 7 &  \winner 6 & -& 15 & 23 &  \winner 5.27505 & 7.54937 & -& 13.85264 & 13.07142 \\
\textsc{ANSYS\_Delor64K\_huber} & 1105173 & 6 & 7 & 7 &  \winner 0 & 3 & 1.30160 & 2.34456 & 2.86475 &  \winner 0.41221 & 1.00721 \\
\textsc{ANSYS\_Delor64K\_lasso} & 7922958 &  \winner 7 &  \winner 7 & -& 17 & 20 & 6.12221 & 9.98326 & -&  \winner 4.06191 & 13.62744 \\
\textsc{HB\_abb313\_huber} & 3748 &  \winner 8 &  \winner 8 & 14 & 11 & 12 &  \winner 0.00124 & 0.00205 & 0.00550 & 0.00456 & 0.00744 \\
\textsc{HB\_abb313\_lasso} & 2887 &  \winner 6 &  \winner 6 & 22 & 8 & 11 &  \winner 0.00085 & 0.00129 & 0.00467 & 0.00434 & 0.00264 \\
\end{longtable}

\end{landscape}

\end{document}